\documentclass[a4paper,reqno]{amsart}
\usepackage{amssymb,epic,eepic} 

\newcommand{\ibar}{\overline{i}}
\newcommand{\jbar}{\overline{j}}
\newcommand{\mvec}[2]{\ensuremath{| #1 \rangle \otimes | #2 \rangle}}
\newcommand{\mvecd}[2]{\ensuremath{\langle #1 | \otimes \langle #2 |}}
\newcommand{\dem}[1]{{\em{#1}}}
\newcommand{\ad}{\text{ad}}
\newcommand{\id}{\text{id}}
\newcommand{\lie}{\mathfrak{g}}
\newcommand{\adjj}{\pi^\Psi}
\newcommand{\adjjc}{\pi^\Psi_{(0)}}
\newcommand{\adj}[1]{\pi^\Psi(#1)}
\newcommand{\adjc}[1]{\pi^\Psi_{(0)}(#1)}
\newcommand{\adjt}[1]{(\pi^\Psi\otimes\pi^\Psi)(\Delta(#1))}
\newcommand{\adjct}[1]{(\pi^\Psi_{(0)}\otimes\pi^\Psi_{(0)})(\Delta(#1))}
\newcommand{\ep}{{\epsilon_1}}
\newcommand{\uqg}{U_h(\mathfrak{g})}
\newcommand{\ug}{U(\mathfrak{g})}
\renewcommand{\a}{\alpha}
\newcommand{\cart}{{\mathcal H}}
\newcommand{\lieh}{\lie_h}
\newcommand{\ch}{\mathbb{C}[[h]]}
\newcommand{\rp}[2]{\langle #1,#2 \rangle}

\theoremstyle{plain}
\newtheorem{proposition}{Proposition}[section]
\newtheorem{lemma}{Lemma}[section]
 
\newtheorem*{dri}{Theorem [Drinfel'd]}
\theoremstyle{definition}
\newtheorem{definition}{Definition}[section]
\numberwithin{equation}{section}
\newcommand{\drawcenteredtext}[3]{\put(#1,#2){\makebox(0,0){#3}}}%
\newcommand{\drawlefttext}[3]{\put(#1,#2){\makebox(0,0)[l]{#3}}}%
\newcommand{\drawrighttext}[3]{\put(#1,#2){\makebox(0,0)[r]{#3}}}%

\newcommand{\drawpath}[4]{\path(#1,#2)(#3,#4)}%

\newcommand{\drawleftbrace}[3]%
{\drawcenteredtext{#1}{#2}{$\left\{ \rule[0mm]{0mm}{#3mm} \right.$}}%

\newcommand{\drawrightbrace}[3]%
{\drawcenteredtext{#1}{#2}{$\left\} \rule[0mm]{0mm}{#3mm} \right.$}}%

\newcommand{\drawoverbrace}[3]%
{\drawcenteredtext{#1}{#2}{$\overbrace{\rule[0mm]{#3mm}{0mm}}$}}%

\newcommand{\drawunderbrace}[3]%
{\drawcenteredtext{#1}{#2}{$\underbrace{\rule[0mm]{#3mm}{0mm}}$}}%

\newcommand{\drawarc}[5]%
{\put(#1,#2){\arc{#3}{#4}{#5}}}%

\begin{document}

\title[The structure of quantum Lie algebras]
{The structure of quantum Lie algebras 
for the classical series B$_l$, C$_l$ and D$_l$}

\author[G.W. Delius]{Gustav W. Delius}
\address{Fakult\"at f\"ur Physik\\
Universit\"at Bielefeld \\
Postfach 10 01 31\\
D-33615 Bielefeld\\
Germany}
\curraddr{Department of Mathematics\\
King's College London\\
Strand\\
London WC2R 2LS\\
Great Britain}
\email{delius@mth.kcl.ac.uk}
\urladdr{http://www.mth.kcl.ac.uk/$\sim$delius/}

\author[C. Gardner]{Christopher Gardner}
\address{Department of Mathematics\\
King's College London\\
Strand\\
London WC2R 2LS\\
Great Britain}
%\curraddr{Institut f\"ur Theoretische Physik\\
%Freie Universit\"at Berlin\\
%Arnimallee 14\\
%D-14195 Berlin\\
%Germany}
\email{cgardner@mth.kcl.ac.uk}
\thanks{C.G. is financially supported by EPSRC}

\author[M.D. Gould]{Mark D. Gould}
\address{Department of Mathematics\\
University of Queensland\\
Brisbane Qld 4072\\
Australia}
\email{mdg@maths.uq.oz.au}

%\dedicatory{Draft dated May 20, 1997}


\begin{abstract}
The structure constants of quantum Lie algebras 
depend on a quantum deformation parameter $q$ and they reduce to 
the classical structure constants of a Lie algebra
at $q=1$. 
We explain the relationship between the structure constants of
quantum Lie algebras and quantum Clebsch-Gordan coefficients for 
adjoint $\otimes$ adjoint $\rightarrow$ adjoint. We present a
practical method for the determination of these quantum Clebsch-Gordan
coefficients and are thus able to give explicit expressions for the
structure constants of the quantum Lie algebras associated 
to the classical Lie algebras B$_l$, C$_l$ and D$_l$. 

In the quantum case also the structure constants of the
Cartan subalgebra are non-zero and we observe that they are
determined in terms of the simple quantum roots.
We introduce an invariant Killing form on the quantum Lie algebras
and find that it takes values which are simple $q$-deformations of
the classical ones. 
\end{abstract}

\maketitle

\section{Introduction\label{sect:intro}}

Quantum Lie algebras are generalizations of Lie algebras whose
structure constants depend on a quantum parameter $q$ and which
are related to the quantized enveloping algebras
(quantum groups) $\uqg$ in a way similar to the way ordinary Lie 
algebras are related to their enveloping algebras $\ug$.

The study of quantum Lie algebras is still in its infancy. There is
no fully developed theory of quantum Lie algebras yet. Instead,
the properties of quantum Lie algebras are being discovered
piecewise through detailed investigation of examples. It is
hoped that these investigations will one day lead to the development
of a full theory. The process is similar to the development of
Lie algebra theory which also began through the detailed study of
the algebras of orthogonal, unitary and symplectic matrices.

In this paper we give the $q$-dependent structure constants of
the quantum Lie algebras associated to the Lie algebras $so(n)$ of special 
orthogonal matrices and the Lie algebras $sp(n)$ of symplectic 
matrices. The case of $sl(n)$ had already been treated in \cite{Del95}.

Previous studies of quantum Lie algebras (cf. \cite{Del96}) have
revealed that one could write the quantum Lie bracket relations
in a form very similar to the classical ones. Namely, choosing
a basis consisting of root generators $X_\alpha$ and Cartan subalgebra
generators $H_i$ the quantum Lie bracket relations are
\begin{eqnarray} \label{qgrad}
\bigl[ H_i , X_{\alpha} \bigr]_h & = 
& l_{\alpha}(H_i)~ X_{\alpha}, \ \ \ ~~\  
\bigl[ X_{\alpha} , H_i \bigr]_h ~~=~~ -r_{\alpha}(H_i)~ X_{\alpha}, 
\nonumber \\
\bigl[H_i , H_j \bigr]_h & = &\sum_k\, f_{ij}{}^{~k}~ H_k, \ \ \ ~ 
\bigl[X_{\alpha} , X_{-\alpha} \bigr]_h ~~=~~-\sum_k\,
g_{\alpha}{}^{~k}~ H_k, \nonumber \\
\bigl[X_{\alpha} , X_{\beta} \bigr]_h & = 
& \left\{ \begin{array}{ll}
N_{\alpha,\beta} ~X_{\alpha+\beta} & 
\mbox{if } \alpha+\beta \text{ is a root} \\
 0 & \mbox{     otherwise}. \end{array} \right. 
\end{eqnarray}
The main differences to the classical relations are:
\begin{itemize}
\item All structure constants now depend on the quantum parameter $q$.
\item Every classical root $\alpha$ splits up into a left quantum
root $l_\alpha$ and a right quantum root $r_\alpha$ which are related
by $q\rightarrow 1/q$.
\item The quantum Lie bracket $[H_i, H_j]_h$
between two Cartan subalgebra generators
is in general non-zero.
\item The quantum Lie bracket is $q$-antisymmetric.
\end{itemize}
As explained in the discussion section, the results in this paper 
allow us to write these quantum Lie bracket
relations using only the $l_\alpha$ and the $N_{\alpha\beta}$.

The notion of a quantum Lie algebra as an ad-submodule of the 
quantized enveloping algebra $U_h(\mathfrak{g})$ which transforms
in the adjoint representation was introduced 
in \cite{Del96}. In that paper it was 
shown that the structure constants of a quantum 
Lie algebra possess a number of symmetries, many of which 
are $q$-generalizations of their classical counterparts. These
symmetries were derived by exploiting the Hopf 
algebra structure of $U_h(\mathfrak{g})$. 

Another approach to quantum 
Lie algebras derives from the notion of bicovariant calculi on quantum 
groups \cite{Wor89,Asc93,Ber90}. There one defines a quantum Lie
product on the dual space 
to the space of left-invariant one-forms, which is an ad-submodule of 
$U_h(\mathfrak{g})$. However in this approach it was not known 
how to ensure that the resulting quantum Lie algebras have the same 
dimension as the corresponding classical Lie algebras except 
in the case of $gl_n$ and, after projection, $A_l$ 
\cite{Fad87,Sud95,DelP96}.
  
In \cite{Del95} the authors constructed the quantum Lie algebras 
associated to $A_l$. They first formed a submodule of $U_h(A_l)$ 
whose elements transform under the adjoint action (defined on 
$U_h(A_l)$) in the vector $\otimes$ dual-vector representation. 
They then projected this module onto the adjoint representation 
inside vector $\otimes$ dual-vector so that the resulting 
ad-module was of the correct dimension. While this approach
can be applied to the construction of quantum Lie algebras associated
to any simple Lie algebra $\lie$, it is a very tedious
way of obtaining the structure constants and we present an
alternative method here.

Our approach to the determination of the structure constants
of quantum Lie algebras relies on the observation \cite{Del97}
that the quantum Lie bracket is an 
intertwiner from adjoint $\otimes$ adjoint $\rightarrow$ adjoint 
where by adjoint we mean the adjoint 
representation of $U_h(\mathfrak{g})$. Thus 
the structure constants are given by the corresponding
inverse quantum Clebsch-Gordan coefficients. In this paper we describe 
a practical procedure for determining these Clebsch-Gordan coefficients.
This allows us to explicitly calculate
the structure constants of the quantum Lie algebras associated 
to the classical Lie algebras B$_l$, C$_l$ and D$_l$.

The outline of this paper is as follows. We begin in section
\ref{sect:qla} by recalling the
definitions and some properties
of quantum Lie algebras. We then explain our method 
for determining their structure constants in section \ref{sect:cgc}. 
In section \ref{sect:kill} we introduce 
an invariant Killing form on the quantum Lie algebras
and observe that on the basis vectors
it takes values which are surprisingly simple deformations of 
the classical ones. Section \ref{sect:results} contains
our results for the structure constants of the quantum Lie 
algebras. By using the Killing form to lower indices we are able
to give concise expressions. 
In section \ref{sect:rels} we derive some further relations between the 
structure constants, which we have verified using Mathematica 
\cite{Wol91} and which provide a good check on the correctness
of our results. In appendix \ref{hopf} we give the basic definition of 
the quantized enveloping algebras $U_h(\mathfrak{g})$ as algebras over
$\mathbb C[[h]]$, the ring of formal power series in $h$.
Further appendices contain a generalized form of Schur's lemma, 
some background material on Clebsch-Gordan coefficients and the 
quantum Lie algebra $(C_2)_h$. 
\vfill\pagebreak[3]

%To fix our notation we write the coproduct  $\Delta$
%on  (associated to the finite-dimensional 
%simple complex Lie algebra $\mathfrak{g}$), whose generators 
%are denoted by $x_i^+,\ x_i^-,\ h_i$, $1 \le i \le 
%\text{rank}(\mathfrak{g})$, as 
%\begin{align} \label{copro}
%\Delta(h_i) &= h_i \hat{\otimes} 1 + 1 \hat{\otimes} h_i, \\
%\Delta(x_i^\pm) &= x_i^\pm \hat{\otimes} q_i^{-h_i/2} +
% q_i^{h_i/2} \hat{\otimes} x_i^\pm.
%\end{align}
%Here $q_i \in \mathbb C[[h]]$, the ring over which 
%$U_h(\mathfrak{g})$ is defined. For further details please consult 
%appendix \ref{hopf}.

\section{Quantum Lie algebras\label{sect:qla}}
\label{qla}
In this section we summarize the definitions for quantum Lie algebras 
given in \cite{Del96,Del95}. In these papers two definitions for 
quantum Lie algebras were given, the resulting algebras being 
isomorphic. The reason for the two definitions stems from the 
fact that we can view a classical Lie algebra $\mathfrak{g}$ 
in two different ways. Either we adopt the perspective that 
$\mathfrak{g}$ is the carrier space for the adjoint 
representation of the classical enveloping algebra
$U(\mathfrak{g})$, this perspective leading to the so called 
abstract quantum Lie algebras $\mathfrak{g}_h$. Or we view 
$\mathfrak{g}$ as a subspace of $U(\mathfrak{g})$.
%on which the Lie bracket is given by the adjoint action. 
As will be seen later, the quantum Lie algebras 
$\mathfrak{L}_h(\mathfrak{g})$ deriving from the second 
definition are just embeddings of the abstract quantum 
Lie algebras $\lie_h$ into $U_h(\mathfrak{g})$.

\subsection{Abstract quantum Lie algebras\label{subsect:aqla}}
\label{qlsec}
Classically a Lie algebra $\mathfrak{g}$ is the carrier space for 
the adjoint representation $\adjjc{}$ of the classical 
enveloping algebra $U(\mathfrak{g})$. 
%, where the classical adjoint 
%representation is defined by $(\ad^{(0)}~ a)~b ~=~ [a,b] ~
%\forall ~a,b~ \in~ \mathfrak{g}$. 
The Lie bracket is a 
$U(\mathfrak{g}$)-module homomorphism from 
$\mathfrak{g}~\otimes\mathfrak{g}$ to $\mathfrak{g}$, i.e.
\begin{equation}
\adjc{x}\circ[~,~] = [~,~]\circ\adjct{x} \ \ \ ~~
\forall x \in U(\mathfrak{g}).
\end{equation}
Furthermore, the following theorem of Drinfel'd \cite{Dri85,Dri90}
provides an isomorphism which allows the construction of the quantum
analog of the adjoint representation.

\begin{dri}
There exists an algebra isomorphism $\varphi~:~ U_h(\mathfrak{g})~ 
\rightarrow U(\mathfrak{g})[[h]]$ such that $\varphi ~\equiv~ 
\id~(mod~h)$ and $\varphi(h_i)~=~h_i$. 
\end{dri}

This implies that if ($V^{\mu},~\pi^{\mu})$ 
denotes a finite-dimensional irreducible representation of 
$U(\mathfrak{g})$ then $(V^{\mu}[[h]],~\pi^{\mu} \circ \varphi)$ 
is a finite-dimensional indecomposable representation of 
$U_h(\mathfrak{g})$. (It is important to recognize that the 
$U_h(\mathfrak{g})$ modules $V[[h]]$ are not irreducible. 
Indeed their submodules are of the form $c~V[[h]]$ with 
$c~ \in~ \mathbb C[[h]]$ not invertible and so Schur's lemma 
is modified, see appendix \ref{schur}.)
In particular $\mathfrak{g}[[h]]$ is a finite-dimensional indecomposable 
module of $U_h(\mathfrak{g})$. Let $\adjj{} = \adjjc{} \circ 
\varphi$ denote the representation of $U_h(\mathfrak{g})$ on 
$\mathfrak{g}[[h]]$. Drinfel'd has shown that this is the only 
way to deform the adjoint representation $\adjjc{}$. With these 
observations in mind, a natural definition for quantum Lie algebras 
is as follows:  

\begin{definition}\label{def:aqla}
A \dem{quantum Lie bracket} is a $U_h(\mathfrak{g}$)-module homomorphism 
$[~,~]_h:\mathfrak{g}[[h]]~\hat{\otimes}~\mathfrak{g}[[h]]~
\rightarrow~\mathfrak{g}[[h]]$ such that $[~,~]_h~=~[~,~]~
(\mbox{mod}~h)$.
\hfill\newline
$\mathfrak{g}_h~=~(\mathfrak{g}[[h]],[~,~]_h)$ is a 
\dem{quantum Lie algebra} (viewed as an algebra over 
$\mathbb C[[h]]$ with non-associative product $[~,~]_h$).
\end{definition}

Such a $U_h(\mathfrak{g})$-module homomorphism $[~,~]_h$
is unique up to scalar 
multiples for $\mathfrak{g} \neq A_l$. The reason for this comes from 
the observation that classically the adjoint representation appears in 
the tensor product of two adjoint representations with unit multiplicity. 
Furthermore, the decomposition of a $U_h(\mathfrak{g})$ tensor product 
representation into indecomposable $U_h(\mathfrak{g})$ modules is 
described by the classical multiplicities of the decomposition of 
the corresponding $U(\mathfrak{g})$ tensor product representation 
into irreducible $U(\mathfrak{g})$ representations. Therefore, 
by the modified form of Schur's lemma, the homomorphism 
$[~,~]_h$ is unique (up to rescaling).

\subsection{Quantum Lie algebras $\mathfrak{L}_h(\mathfrak{g})$ 
inside $U_h(\mathfrak{g})$\label{subsect:qlai}}

Alternatively $\mathfrak{g}$ can be viewed as a subspace of its 
enveloping algebra $U(\mathfrak{g})$ with the Lie bracket on this 
subspace given by the adjoint action of $U(\mathfrak{g})$. So 
another natural definition of a quantum Lie algebra is as 
an ad-submodule of $U_h(\mathfrak{g})$ with the quantum Lie bracket 
given by the adjoint action of $U_h(\mathfrak{g})$.

\begin{definition}\label{def2}
A \dem{quantum Lie algebra} $\mathfrak{L}_h(\mathfrak{g})$ 
\dem{inside} $U_h(\mathfrak{g})$ is a finite-di\-men\-sio\-nal 
indecomposable ad-submodule of $U_h(\mathfrak{g})$ endowed 
with a \dem{quantum Lie bracket} $[a , b]_h~=~(\ad~a)~b$ 
such that: \\
(a) $\mathfrak{L}_h(\mathfrak{g})$ is a deformation of $\mathfrak{g}$, 
i.e., there is an algebra isomorphism 
$\mathfrak{g}\cong$ $\mathfrak{L}_h(\mathfrak{g})|_{h=0}$, \\
(b) $\mathfrak{L}_h(\mathfrak{g})$ is invariant under 
$\tilde{\theta},\ \tilde{S}$ and any diagram automorphism $\tau$,
where $\sim$ denotes $q$-conjugation, $\theta$ the Cartan involution 
and $S$ the antipode, see appendix \ref{hopf}. \\
A \dem{weak quantum Lie algebra} $\mathfrak{l}_h(\mathfrak{g})$ 
is defined similarly but without the second condition.
\end{definition}

By ad-submodule we mean that it is invariant under the adjoint action 
of $U_h(\mathfrak{g})$ which, in Sweedler's notation \cite{Swe69}, 
is given by
\begin{equation} 
(\ad~x)~y ~~=~~ \sum~x_{(1)} ~y~S(x_{(2)}),\ \ \ ~~~x,y~\in U_h(\mathfrak{g}).
\end{equation}
This adjoint action defines an infinite dimensional representation 
of $U_h(\mathfrak{g})$ on itself.  Notice that classically 
this definition for the adjoint action reduces to the usual commutator 
when restricted to the Lie algebra naturally embedded in 
$U(\mathfrak{g})$, giving rise to the classical adjoint representation 
$\adjjc$.

\begin{proposition} \label{prop1} The adjoint action restricted to 
$\mathfrak{L}_h(\mathfrak{g}) \otimes \mathfrak{L}_h(\mathfrak{g})$ 
is a $U_h(\mathfrak{g})$-module homomorphism from 
$\mathfrak{L}_h(\mathfrak{g}) \otimes \mathfrak{L}_h(\mathfrak{g})$ 
to $\mathfrak{L}_h(\mathfrak{g})$.
\end{proposition} 
\begin{proof}
We need to show that
\begin{equation} \label{prfq}
\sum \biggl( (\ad(\ad~x_{(1)})~a)((\ad~x_{(2)})~(b)) \biggr) =
(\ad~x)((\ad~a)~(b)),\ \ \ \forall a,b\in \mathfrak{L}_h(\mathfrak{g})
\nonumber
\end{equation}
which, using that $[a , b]_h~=~(\ad~a)~b$, is equivalent to
\begin{equation}
(\ad~x) ~\circ~[ ~,~ ]_h=~~ [~,~]_h~\circ~ 
(\ad ~ \otimes ~\ad ~ )(\Delta(x)),\ \ \ ~~~ \forall x ~ \in~ 
U_h(\mathfrak{g}).
\nonumber
\end{equation}
The left hand side of eq. \eqref{prfq}, using cocommutativity of the 
Hopf algebra, can be expressed as
\begin{eqnarray*}
&&\sum (\ad x_{(1)}aS(x_{(2)}))(x_{(3)}bS(x_{(4)})) \\
&&~~~~~ = \sum x_{(1)}a_{(1)}S(x_{(4)})x_{(5)}bS(x_{(6)})
S(x_{(2)}a_{(2)}S(x_{(3)})) \\
&&~~~~~ =  \sum x_{(1)}a_{(1)} \epsilon(x_{(4)})bS(x_{(5)})
S^2(x_{(3)})S(a_{(2)})S(x_{(2)}) \\
&&~~~~~ =  \sum x_{(1)}a_{(1)}bS(a_{(2)})S(x_{(2)}) \\
&&~~~~~ =  (\ad x)((\ad a)b).
\end{eqnarray*}
where we have used the Hopf algebra property  
$S(y_{(1)})y_{(2)} = \epsilon(y)$, and that $S$ is a Hopf algebra 
anti-automorphism.
\end{proof}

How are these quantum Lie algebras related to the abstract quantum 
Lie algebras defined in the previous section? This question is 
answered in the following proposition.

\begin{proposition} All weak quantum Lie algebras 
$\mathfrak{l}_h(\mathfrak{g})$ inside $U_h(\mathfrak{g})$ 
are isomorphic to an abstract quantum Lie algebra $\mathfrak{g}_h$.
\end{proposition}
\begin{proof}
$\mathfrak{l}_h(\mathfrak{g})$ are finite-dimensional 
indecomposable $U_h(\mathfrak{g})$ modules. By condition (a) of
definition \ref{def2}, 
they carry a deformation of the representation of $U(\mathfrak{g})$ 
carried by $\mathfrak{g}$. As mentioned earlier, Drinfel'd showed 
that there is only one such deformation of the classical adjoint 
representation $\adjjc{}$, namely the adjoint representation 
$\adjj{}=\adjjc\circ\varphi$ carried by $\mathfrak{g}[[h]]$. Therefore 
$\mathfrak{l}_h(\mathfrak{g})$ is isomorphic to $\mathfrak{g}[[h]]$ 
as a $U_h(\mathfrak{g})$-module. Furthermore, by proposition 
\eqref{prop1} the product on $\mathfrak{l}_h(\mathfrak{g})$ is a 
$U_h(\mathfrak{g})$-module homomorphism.
\end{proof}
In other words, the weak quantum Lie algebras  
$\mathfrak{l}_h(\mathfrak{g})$ are those embeddings of 
$\mathfrak{g}[[h]]$ in $U_h(\mathfrak{g})$ on which the adjoint 
action of $U_h(\mathfrak{g})$ coincides with the adjoint 
representation $\adjj{}$.

\subsection{General properties of quantum Lie algebras}

\begin{proposition} Given the grading of the classical Lie algebra 
by the set of non-zero roots $R$ and zero, described by 
\begin{eqnarray}
\mathfrak{g} & = & \bigoplus_{\alpha \in R}~\mathfrak{g}_{\alpha} ~
\oplus~ \mathfrak{g}_0, \ \ \ \ \ \ [\mathfrak{g}_{\alpha},
\mathfrak{g}_{\beta}]\subset \mathfrak{g}_{\alpha+\beta}
\end{eqnarray}
with $\mathfrak{g}_{\alpha}~=~ \{x \in \mathfrak{g} \mid 
\adjc{h_i}x~=~\alpha(h_i)x ~\forall h_i \}$, a quantum 
Lie algebra $\mathfrak{g}_h$ possesses the grading
\begin{eqnarray}
\mathfrak{g}_h & = & \bigoplus_{\alpha~ \in~ R}~
\mathfrak{g}_{\alpha}[[h]] ~\oplus~ \mathfrak{g}_o[[h]],\ \ \  
\bigl[\mathfrak{g}_{\alpha}[[h]] , \mathfrak{g}_{\beta}[[h]] 
\bigr]_h \subset \mathfrak{g}_{\alpha+\beta}[[h]]. 
\end{eqnarray}
\end{proposition}

\begin{proof}
The isomorphism $\varphi$ leaves $h_i$ invariant. 
This implies $\adj{h_i} = \adjc{\varphi(h_i)} = \adjc{h_i}$, so that
\begin{equation}
\mathfrak{g}_{\alpha}[[h]] ~~=~~ \{ x~ \in \mathfrak{g}[[h]] 
\mid \adj{h_i}~x ~=~ \alpha(h_i)~x ~\forall ~h_i \}. 
\end{equation}
%Indeed $\mathfrak{g}_{\alpha}[[h]] = \mathfrak{g}_{\alpha} \otimes 
%\mathbb C[[h]]$ because if $\mathfrak{g}_{\alpha}[[h]]$ contains 
%an element $y \in \mathfrak{g}_{\beta}, ~~\beta \neq \alpha$ 
%then $\adj{h_i}~y ~=~ \beta(h_i) y$ in other words 
%$y \in \mathfrak{g}_{\beta}[[h]]$.
Let $X_{\alpha} \in \mathfrak{g}_{\alpha}[[h]]$ and $X_{\beta} \in 
\mathfrak{g}_{\beta}[[h]]$. By definition $[~,~]_h$ is a 
$U_h(\mathfrak{g})$-module homomorphism and so
\begin{eqnarray}
\adj{h_i}~ [X_{\alpha} , X_{\beta}]_h & =
& [~,~]_h \circ \adjt{h_i}~ (X_{\alpha} \hat{\otimes} 
X_{\beta})  \nonumber \\
& = & [\adj{h_i}~X_{\alpha} , X_{\beta}]_h ~+~ 
[X_{\alpha} , \adj{h_i}~X_{\beta}]_h \nonumber \\
& = & (\alpha(h_i)+\beta(h_i)) [X_{\alpha} , X_{\beta}]_h,
\end{eqnarray}
thus $[X_{\alpha} , X_{\beta}]_h \in 
\mathfrak{g}_{\alpha+\beta}[[h]]$. 
%Here $\ad_2^{(h)}(h_i) ~=~ (\ad^{(h)} \otimes \ad^{(h)}) \Delta(h_i)$.
\end{proof}

If we now choose a basis for the quantum Lie algebra $\mathfrak{g}_h$, 
given by $\{X_{\alpha} \in \mathfrak{g}_{\alpha} \mid 
\alpha \in R \} \cup \{H_i  \in \mathfrak{g}_0 \mid i=1,\dots
\mbox{rank($\mathfrak{g}$)} \}$, then, because of the grading, 
the Lie bracket is restricted to the following form:
\begin{eqnarray}
\bigl[ H_i , X_{\alpha} \bigr]_h & = 
& l_{\alpha}(H_i)~ X_{\alpha}, \ \ \ \  
\bigl[ X_{\alpha} , H_i \bigr]_h ~~=~~ -r_{\alpha}(H_i)~ X_{\alpha}, 
\nonumber \\
\bigl[H_i , H_j \bigr]_h & = &\sum_k\, f_{ij}^{~~k}~ H_k, \ \ \ ~
\bigl[X_{\alpha} , X_{-\alpha} \bigr]_h ~~=~~-\sum_k\,
g_{\alpha}^{~k}~ H_k, \nonumber \\
\bigl[X_{\alpha} , X_{\beta} \bigr]_h & = 
& \left\{ \begin{array}{ll}
N_{\alpha,\beta} ~X_{\alpha+\beta} & 
\mbox{    for  } \alpha+\beta ~\in~ R \\
 0 & \mbox{     otherwise}. \end{array} \right. 
\end{eqnarray}
We call the zero-weight subalgebra $\cart:=\lie_0[[h]]$ of 
$\lieh$ spanned by the 
generators $\{H_i\}_{i=1,\dots,\text{rank}\lie}$ the Cartan subalgebra
even though unlike in the classical case $[H_i,H_j]_h\neq 0$. 
%We also denote this Cartan subalgebra by $\cart$.
The quantum roots $l_\alpha$ and $r_\alpha$ are linear functionals
on $\cart$. The structure constants $l_\alpha(H_i)$,
$r_\alpha(H_i)$, $g_\alpha^k$, 
$f_{ij}{}^k$ and $N_{\alpha\beta}$, are power series in the
indeterminate $h$ (in fact most of them turn out to be polynomials
in $q^{1/2}=e^{h/2}$).

The quantum structure constants possess a number of symmetries 
\cite{Del96,Del97}, many of which are natural $q$-generalizations 
of their classical counterparts. The basis can be chosen so that 
the symmetries are given by:
\begin{gather}
l_{\alpha} =  \tilde{r}_{\alpha},\ \ \ \ \ \ ~
f_{ij}{}^k=-\tilde{f}_{ji}{}^k,\nonumber \\ 
N_{\alpha,\beta}=-\tilde{N}_{\beta,\alpha},\ \ \ \ \ \ ~
g_{\alpha}^{~k}=-\tilde{g}_{-\alpha}^{~k},  
\label{qant1}
\end{gather}
\begin{gather}
l_{\alpha}=-\tilde{l}_{-\alpha},\ \ \ \ \ \ ~
f_{ij}^{~~k}=-\tilde{f}_{ij}^{~~k},\nonumber \\
N_{\alpha,\beta}=-\tilde{N}_{-\alpha,-\beta}, 
\label{qant2}
\end{gather}
\begin{align}
N_{\alpha,\beta} &= q^{\rho \cdot \beta} N_{\beta,-\alpha-\beta},
\label{qant3} \\
\sum_l f_{jk}^{~~l} B(H_i,H_l) & =  \sum_l f_{ji}^{~~l} B(H_l,H_k),
\label{qant4}\\
-g_{\alpha}^{~i} B(H_i,H_j) & = 
l_{\alpha}(H_j) q^{-\rho\cdot\alpha}, 
\label{qant5}
\end{align}
where $B$ is the Killing form, to be discussed in section
\ref{sect:kill}, $\rho$ is half 
the sum of the positive roots, and $\sim$ denotes $q$-conjugation,
%$\mbox{$\sim: \mathbb C[[h]] \rightarrow \mathbb C[[h]]$}, 
%$a\mapsto\tilde{a}$ is the ring automorphism defined by 
%$\tilde{h}=-h$. 
i.e., the operation $h\mapsto -h$ or equivalently $q\mapsto 1/q$.
Note that in this paper we choose the basis for the quantum Lie algebra
slightly differently from \cite{Del96} to satisfy
\begin{gather}
B(X_\alpha,X_{-\alpha})=q^{-\rho\cdot\alpha},\\
\tilde{\theta}(X_{\alpha})=-X_{-\alpha},\ \ \ \ \ \ \ \ \ ~
\tilde{\theta}(H_i)=-H_i,\\
\tilde{S}(X_\alpha)=-q^{\rho\cdot\alpha}\,X_\alpha,\ \ \ \ \ \ \ \ \ ~
\tilde{S}(H_i)=-H_i.
\end{gather}

The symmetries in \eqref{qant1} 
express the $q$-antisymmetry of the quantum Lie bracket. 
\begin{equation}
[a , b ]_h^{\triangledown} = 
 - [b^{\triangledown} , a^{\triangledown}]_h.
\end{equation}
The map $^{\triangledown}$:$\mathfrak{g}[[h]] \rightarrow \mathfrak{g}[[h]]$ 
is defined so that, if $a=\sum_{\alpha} \zeta_{\alpha} X_{\alpha} + 
\sum_i \eta_i H_i$ is a general element of $\mathfrak{g}[[h]]$, 
then $a^{\triangledown} = \sum_{\alpha} \tilde{\zeta}_{\alpha} X_{\alpha} + 
\sum_i \tilde{\eta}_i H_i$.
In the classical limit $q$-antisymmetry expresses the antisymmetry 
of the Lie bracket. 
The origin of the $q$-antisymmetry of the 
quantum Lie bracket was explained in \cite{Del97}.

\subsection{Quantum structure constants as quantum
Clebsch-Gordan coefficients\label{subsect:qsc}}

Let $\{v^a\}_{a=1,\dots\text{dim}\lie}$ be a basis for $\mathfrak{g}[[h]]$. 
The highest weight vector $v^1 $ generating this 
$U_h(\mathfrak{g})$-module satisfies the relations,
\begin{equation} \label{adrel}
\adj{x_i^+}~v^1  ~=~0,\ \ \ \ \ \ \adj{h_i}~v^1 ~=~\Psi(h_i)~v^1 ,\ \ \  
\forall~i,
\end{equation}
where $\Psi$ is the highest root of $\mathfrak{g}$ and 
$\adjj{}~=~\adjjc{} \circ \varphi$. Let $P^a(x^-)$ be polynomials in 
the $x_i^-$ such that $v^a =\adj{P^a(x^-)}~v^1 $. The adjoint 
representation matrices in this basis are given by
\begin{equation} \label{repm}
\adj{x}~v^a ~~=~~v^b ~\pi_b^{\Psi a}(x).
\end{equation}
Notation: latin indices which appear as both upper and lower indices
in the same expression are summed over.

We know that $\mathfrak{g}[[h]]~\hat{\otimes}~\mathfrak{g}[[h]]$ 
contains $\mathfrak{g}[[h]]$ with unit multiplicity for the 
algebras $B_l$, C$_l$ and $D_l$. Hence, there exists a 
highest weight state $\hat{v}^1$ inside 
$\mathfrak{g}[[h]]~\hat{\otimes}~\mathfrak{g}[[h]]$ that 
satisfies the analogue of relations \eqref{adrel}, namely
\begin{equation}
\adjt{x_i^+}~\hat{v}^1 ~=~0,\ \ \ ~
\adjt{h_i}~\hat{v}^1 ~=~\Psi(h_i)~\hat{v}^1 ,~~\forall~i.
\end{equation}
and $\hat{v}^1$ is the unique (up to rescaling) state with this property.
%We construct a basis $\{\hat{v}^a \}$ with the same polynomials 
%$P^a(x^-)$. 
The vectors $\hat{v}^a =\adjt{P^a(x^-)}~\hat{v}^1$
form a basis for 
$\mathfrak{g}[[h]]$ inside $\mathfrak{g}[[h]]~\hat{\otimes}~
\mathfrak{g}[[h]]$ such that
\begin{equation}
\adjt{x}~\hat{v}^a ~~=~~\hat{v}^b ~\pi_b^{\Psi a}(x)
\end{equation}
with the same representation matrices as in \eqref{repm}.
We can expand the vectors $\hat{v}^a$ in terms of the tensor
product basis as
\begin{equation}
\hat{v}^a=~v^b\otimes v^c~C^{\Psi}_b{}^\Psi_c |^a_\Psi.
\end{equation}
The $C^{\Psi}_b{}^\Psi_c |^a_\Psi\in \mathbb C[[h]]$ are the 
Clebsch-Gordan coefficients. Notation: we often indicate above or
below an index which representation the index belongs to.
We do this by giving
the heighest weight of the representation, e.g., in 
$C^{\Psi}_b{}^\Psi_c |^a_\Psi$ all indices belong to the
representation with highest weight $\Psi$ ($\Psi$ being the highest
root of the Lie algebra), i.e., the adjoint representation.

The embedding map $\beta:~\mathfrak{g}[[h]] \rightarrow 
\mathfrak{g}[[h]]~\hat{\otimes}~\mathfrak{g}[[h]]$, defined by 
$v^a  \rightarrow \hat{v}^a $ is a $U_h(\mathfrak{g})$-module 
homomorphism, i.e. 
$\beta ~\circ~ \adj{x}~=~\adjt{x} ~\circ~ \beta$.
This is easy to check:
\begin{eqnarray*}
(\beta ~\circ~ \adj{x})~v^a  & 
= & \beta(v^b ~\pi_b^{\Psi a}(x))~~=~~ \hat{v}^b  ~\pi_b^{\Psi a}(x) \\
& = & \adjt{x}~\hat{v}^a ~~=~~(\adjt{x}~\circ~\beta)~v^a .
\end{eqnarray*}
Both $\mathfrak{g}[[h]]$ and Im$~\beta$ are indecomposable modules 
and so $\beta$ (with its range restricted to Im$~\beta$) is unique 
and invertible by the weak form of Schur's lemma. We now define the 
quantum Lie bracket $[~,~]_h:~\mathfrak{g}[[h]] ~\hat{\otimes}~
\mathfrak{g}[[h]] \rightarrow \mathfrak{g}[[h]]$ to be zero on the 
module complement of Im$~\beta$ whilst on Im$~\beta$ we define 
$[~,~]_h ~=~ \beta^{-1}$. Then
\begin{equation}
[v^a  , v^b ]~~=~~C^\Psi_{c} |^a_\Psi{}^b_\Psi~v^c , \ \ \ ~
\mbox{with } C^{\Psi}_{d} |^b_\Psi{}^c_\Psi~
C_b^\Psi{}_c^\Psi|^a_\Psi ~=~\delta_d^a.
\end{equation}
Finally, because $\beta$ is a $U_h(\mathfrak{g})$-module homomorphism, 
$\beta^{-1}$ is also a $U_h(\mathfrak{g})$-module homomorphism. So we 
have shown that we can construct a bracket with the correct properties.

\dem{Conclusion}: to construct the quantum Lie algebra $\mathfrak{g}_h$ 
associated to $\mathfrak{g}$ we need to calculate the inverse 
Clebsch-Gordan coefficients for the decomposition 
$\mathfrak{g}[[h]]~\hat{\otimes}~\mathfrak{g}[[h]]$ into 
$\mathfrak{g}[[h]]$.

\section{\sloppy{Calculation of quantum Clebsch-Gordan co\-efficients}
\label{sect:cgc}}
In this section we present our method for the calculation of
the inverse quantum
Clebsch-Gordan coefficients for adjoint $\otimes$ adjoint 
$\rightarrow$ adjoint. We don't actually form adjoint $\otimes$ 
adjoint. Instead we build the Clebsch-Gordan coefficients from 
vector $\otimes$ vector. The reason for this lies in the 
relative simplicity of the vector representations of 
$U_h(\mathfrak{g})$ for $B_l$, $C_l$ and $D_l$. The most attractive 
feature of the vector representation is that the classical 
representation matrices for the Cartan subalgebra generators 
have eigenvalues $+1,0\mbox{ or }-1$. This means that the 
classical representation matrices are also representation 
matrices for $U_h(\mathfrak{g})$.

Let us first remind ourselves of the roots for the algebras 
$B_l$, $C_l$ and $D_l$. 
Let $\epsilon_i$ denote the orthonormal basis vectors of the root
space $\Phi$.
The roots for $B_l$ can be written as 
$\pm(\epsilon_i \pm \epsilon_j)$ with $i \neq j$ and 
$\pm \epsilon_i$ with $i \le l$. 
The simple roots are $\alpha_i = \epsilon_i-\epsilon_{i+1}$ for 
$i < l$ and $\alpha_l = \epsilon_l$. 
For $C_l$ the roots are given by $2 \epsilon_i$ with $i \le l$ and 
$\pm(\epsilon_i \pm \epsilon_j)$ with $i \neq j$. The simple roots 
are $\alpha_i=\epsilon_i-\epsilon_{i+1}$ for $i < l$ and 
$\alpha_l = 2 \epsilon_l$. 
For $D_l$ the roots are given by $\pm(\epsilon_i \pm \epsilon_j)$ 
with $i \neq j$ and the simple roots are 
$\alpha_i=\epsilon_i-\epsilon_{i+1}$ for $i < l$ and 
$\alpha_l = \epsilon_{l-1}+ \epsilon_l$. 
The vector representation is generated from the highest weight state, 
denoted by $| 1 \rangle$, whose weight is $\epsilon_1$. 
This corresponds to $\alpha_1+\alpha_2+...+\alpha_l$ for 
$B_l$, $\alpha_1+\alpha_2+...+\frac{1}{2} \alpha_l$ for $C_l$, and 
$\alpha_1+\alpha_2+...+\frac{1}{2}(\alpha_{l-1}+\alpha_l)$ for $D_l$. 
The value of $\rho \cdot \alpha$ is calculated easily from the 
value of $\rho \cdot \epsilon_i$ ($i \leq l$) as given below.
\begin{equation}
\rho \cdot \epsilon_i ~=~ \left\{ \begin{array}{ll}
l-i+1/2 & \mbox{    for $B_l$ } \\
l-i+1 & \mbox{    for $C_l$ } \\
l-i & \mbox{    for $D_l$ } 
\end{array} \right.
\end{equation}

In the following we will denote by $V^\mu$ the representation space 
with highest weight $\mu$, e.g., $V^\Psi$ is the adjoint representation,
$V^\ep$ the vector representation and $V^0$ the singlet representation.

\subsection{The Clebsch-Gordan coefficients for adjoint 
$\otimes$ adjoint $\rightarrow$ adjoint\label{subsect:cgca}}
\label{secrts}

To calculate these Clebsch-Gordan coefficients we start from the 
simpler decomposition of the 
direct product representation $V^\ep\otimes V^\ep
\rightarrow V^{2 \ep}\oplus V^{\ep+\epsilon_2} \oplus V^0$.
We denote by $\{v^a_\mu\}_{a=1,\dots,\text{dim}V^\mu}$
the basis of a representation space $V^\mu$. Then the decomposition 
is described by the inverse Clebsch-Gordan coefficients as follows
\begin{equation}\label{exp}
v^a_{\epsilon_1}\otimes v^b_{\epsilon_1}=
\sum_{c=1}^{\text{dim}V^{2\epsilon_1}}\ v_{2\epsilon_1}^c\ 
C_c^{2\ep}|^a_\ep{}^b_\ep
+\sum_{c=1}^{\text{dim}V^{\ep+\epsilon_2}}\ v_{\ep+\epsilon_2}^c\ 
C_c^{\ep+\epsilon_2}|^a_\ep{}^b_\ep
+v_0\ C^0|^a_\ep{}^b_\ep.
\end{equation}

It is very convenient to introduce a graphical notation for
the Clebsch-Gordan coefficients and their inverses as in 
figure \ref{fig1}. 
\begin{figure} 
% Drawing generated by LaTeX-CAD 1.7 - requires latexcad.sty 
% (c) 1996 John Leis leis@usq.edu.au 
\begin{picture}(145,78)
\thinlines
\drawarc{60.0}{44.0}{20.0}{3.14}{6.28}
\drawpath{60.0}{54.0}{60.0}{64.0}
\drawrighttext{58.0}{60.0}{$\lambda$}
\drawrighttext{50.0}{50.0}{$\mu$}
\drawlefttext{70.0}{50.0}{$\nu$}
\drawcenteredtext{60.0}{68.0}{$c$}
\drawcenteredtext{50.0}{40.0}{$a$}
\drawcenteredtext{70.0}{40.0}{$b$}
\drawrighttext{46.0}{54.0}{$C_a^\mu{}_b^\nu|^c_{\lambda}\quad=\quad$}
\drawrighttext{118.0}{54.0}{$C_c^{\lambda}|^a_\mu{}^b_\nu\quad=\quad$}
\drawarc{132.0}{64.0}{20.0}{6.28}{3.14}
\drawpath{132.0}{54.0}{132.0}{44.0}
\drawcenteredtext{122.0}{68.0}{$a$}
\drawcenteredtext{142.0}{68.0}{$b$}
\drawcenteredtext{132.0}{40.0}{$c$}
\drawlefttext{142.0}{58.0}{$\nu$}
\drawrighttext{122.0}{58.0}{$\mu$}
\drawrighttext{130.0}{48.0}{$\lambda$}
\drawrighttext{46.0}{20.0}{$C_a^\mu{}_b^\mu|_0\quad=\quad$}
\drawarc{60.0}{14.0}{20.0}{3.14}{6.28}
\drawrighttext{118.0}{20.0}{$C^0|^a_\mu{}^b_\mu\quad=\quad$}
\drawarc{132.0}{24.0}{20.0}{6.28}{3.14}
\drawcenteredtext{50.0}{10.0}{$a$}
\drawcenteredtext{70.0}{10.0}{$b$}
\drawcenteredtext{142.0}{28.0}{$b$}
\drawcenteredtext{122.0}{28.0}{$a$}
\drawrighttext{52.0}{22.0}{$\mu$}
\drawlefttext{68.0}{22.0}{$\mu$}
\drawrighttext{124.0}{16.0}{$\mu$}
\drawlefttext{140.0}{16.0}{$\mu$}
\end{picture}
\caption{Graphical representation of Clebsch-Gordan coefficients
\label{fig1}}
\end{figure}

The composition of any number of 
intertwiners is again an intertwiner. We exploit this
fact to build the intertwiner for adjoint $\otimes$ adjoint 
$\rightarrow$ adjoint from the intertwiners for
vector $\otimes$ vector into adjoint and singlet as shown
in the following diagram
\begin{equation}
% Drawing generated by LaTeX-CAD 1.7 - requires latexcad.sty 
% (c) 1996 John Leis leis@usq.edu.au 
\begin{picture}(84,44)
\thinlines
\drawarc{20.0}{28.0}{16.0}{6.28}{3.14}
\drawpath{52.0}{34.0}{52.0}{30.0}
\drawarc{60.0}{26.0}{24.0}{6.28}{3.14}
\drawarc{52.0}{26.0}{8.0}{3.14}{6.28}
\drawarc{60.0}{26.0}{8.0}{6.28}{3.14}
\drawarc{68.0}{26.0}{8.0}{3.14}{6.28}
\drawpath{68.0}{34.0}{68.0}{30.0}
\drawrighttext{58.0}{10.0}{$\Psi$}
\drawlefttext{74.0}{26.0}{$\ep$}
\drawpath{20.0}{20.0}{20.0}{12.0}
\drawrighttext{46.0}{26.0}{$\ep$}
\drawlefttext{70.0}{32.0}{$\Psi$}
\drawrighttext{10.0}{28.0}{$\Psi$}
\drawcenteredtext{12.0}{32.0}{$a$}
\drawrighttext{50.0}{32.0}{$\Psi$}
\drawcenteredtext{28.0}{32.0}{$b$}
\drawcenteredtext{52.0}{38.0}{$a$}
\drawcenteredtext{68.0}{38.0}{$b$}
\drawlefttext{30.0}{28.0}{$\Psi$}
\drawcenteredtext{60.0}{26.0}{$\ep$}
\drawpath{60.0}{14.0}{60.0}{10.0}
\drawpath{48.0}{32.0}{48.0}{32.0}
\drawpath{48.0}{32.0}{48.0}{32.0}
\drawrighttext{18.0}{16.0}{$\Psi$}
\drawcenteredtext{20.0}{8.0}{$c$}
\drawcenteredtext{60.0}{6.0}{$c$}
\drawcenteredtext{38.0}{20.0}{$=$}
\end{picture}
\end{equation}
or equivalently
\begin{equation}
\label{intw1}
C_c^\Psi|^a_\Psi{}^b_\Psi\ =\ C_i^\ep{}_j^\ep|^a_\Psi\ 
C_k^\ep{}_l^\ep|^b_\Psi\ C^0|^j_\ep{}^k_\ep\ C_c^\Psi|^i_\ep{}^l_\ep.
\end{equation}
Let us show explicitly that the $C^\Psi|_{\Psi\Psi}$ defined
by this equation does indeed satisfy the intertwining property
\begin{equation} \label{intw}
C_c^\Psi|^{a'}_\Psi{}^{b'}_\Psi\ \pi^\Psi_{a'}{}^a(x_{(1)})\ 
\pi^\Psi_{b'}{}^b(x_{(2)})\ =\ 
\pi^\Psi_c{}^{c'}(x)\ C^\Psi_{c'}|^a_\Psi{}^b_\Psi \ \ \ ~\forall x
\in U_h(\mathfrak{g}).
\end{equation}
Substituting \eqref{intw1} and using the intertwining property of 
$C^{\ep\ep}|_\Psi$, the left hand side of \eqref{intw} 
becomes
\begin{equation}
\pi^\ep_i{}^{i'}(x_{(1)})\ \pi^\ep_j{}^{j'}(x_{(2)})\ 
C_{i'}^\ep{}_{j'}^\ep|^a_\Psi\ 
\pi^\ep_k{}^{k'}(x_{(3)})\ \pi^\ep_l{}^{l'}(x_{(4)})\ 
C_{k'}^\ep{}_{l'}^\ep|^b_\Psi\ 
C^0|^j_\ep{}^k_\ep\ C_c^\Psi|^i_\ep{}^l_\ep.
\end{equation}
The 1-dimensional representation of $U_h(\mathfrak{g})$ is given by 
the counit $\epsilon$ (recall $\epsilon(xy)=\epsilon(x)\epsilon(y)$ 
$~\forall~x,y\in U_h(\mathfrak{g})$), and so
\begin{equation}
C^0|^j_\ep{}^k_\ep\ \pi^\ep_j{}^{j'}(x_{(2)})\ \pi^\ep_k{}^{k'}(x_{(3)})
\ =\ \epsilon(x_{(2)})\ C^0|^{j'}_\ep{}^{k'}_\ep
\end{equation}
Consequently one obtains (using cocommutativity)
\begin{equation}
C_{i'}^\ep{}_{j'}^\ep|^a_\Psi\ C_{k'}^\ep{}_{l'}^\ep|^b_\Psi\ 
C^0|^{j'}_\ep{}^{k'}_\ep\ C_c^\Psi|^i_\ep{}^l_\ep\ 
\pi^\ep_i{}^{i'}(x_{(1)})\ \epsilon(x_{(2)})\ 
\pi^\ep_l{}^{l'}(x_{(3)})
\end{equation}
which is equivalent to the right hand side of \eqref{intw}, 
if we use the intertwining property of the Clebsch-Gordan 
coefficient $C^\Psi|_{\ep\ep}$ and the counit property 
$\epsilon(x_{(1)})x_{(2)}=x$.

\subsection{The vector representation $V^\ep$\label{subsect:vector}}

\label{secvec}
We now turn our attention to the vector representation. 
We change notation and denote the basis states for the vector 
representation by $| i\rangle$ with $i=1,...,2l$ for 
$C_l$ and $D_l$, and $i=1,...,2l+1$ for $B_l$. The matrices 
$e_{ij}$ act on these basis states as $e_{ij} | k \rangle = 
\delta_{jk} | i \rangle$.

The states in the vector representation have the following weights:
\begin{eqnarray*}
| i \rangle & \longleftrightarrow & \epsilon_i, \ \ \ \ \ \ i \leq l, \\
| l+1 \rangle & \longleftrightarrow & 0, \\
| i \rangle & \longleftrightarrow & -\epsilon_{\ibar}, \ \ \ \ \ \ 
l+1 < i \leq 2l+1,
\end{eqnarray*}
for the algebra $B_l$ ($\ibar=2l+2-i$), and
\begin{eqnarray*}
| i \rangle & \longleftrightarrow & \epsilon_i, \ \ \ \ \ \ i \leq l, \\
| i \rangle & \longleftrightarrow & -\epsilon_{\ibar}, \ \ \ \ \ \ 
l < i \leq 2l,
\end{eqnarray*}
for the algebras $C_l$ and $D_l$ ($\ibar = 2l+1-i$).

In the 
following the $E_i$, $F_i$ and $H_i$ are short hand for 
$\pi^{\epsilon_1}(x_i^{+})$, $\pi^{\epsilon_1}(x_i^{-})$ and 
$\pi^{\epsilon_1}(h_i)$. 
\subsubsection{$\mathfrak{g}~=~B_{l}$}
\begin{eqnarray}
&&E_i=a_{i,i+1},\ \ \ ~
F_i=a_{i+1,i},\ \ \ ~
H_i=a_{ii}-a_{i+1,i+1},\ \ \ ~1\leq i < l,
 \nonumber \\
&&E_l=\sqrt{2}a_{l,l+1},\ \ \ ~
F_l=\sqrt{2}a_{l+1,l},\ \ \ ~
H_l=2a_{l,l},
\end{eqnarray}
where
\begin{equation}
a_{ij} = e_{ij}-e_{\jbar\ibar},\ \ \ ~
\ibar=2l+2-i. 
\end{equation}

\subsubsection{$\mathfrak{g}~=~C_{l}$}
\begin{eqnarray}
&&E_i=a_{i,i+1},\ \ \ ~
F_i=a_{i+1,i},\ \ \ ~
H_i=a_{ii}-a_{i+1,i+1},\ \ \ ~1\leq i < l,
\nonumber \\
&&E_l=\frac{1}{2}a_{l,l+1},\ \ \ ~
F_l=\frac{1}{2}a_{l+1,l},\ \ \ ~
H_l=a_{l,l},
\end{eqnarray} 
where
\begin{eqnarray*}
a_{ij}=e_{ij}-(-1)^{i+j}e_{\jbar\ibar},\ \ \ ~
\ibar=2l+1-i. 
\end{eqnarray*}

\subsubsection{$\mathfrak{g}~=~D_{\mathit l}$} 
\begin{eqnarray}
&&E_i=a_{i,i+1},\ \ \ ~
F_i=a_{i+1,i},\ \ \ ~
H_i=a_{ii}-a_{i+1,i+1},\ \ \ ~1\leq i < l,
\nonumber\\
&&E_l=a_{l-1,l+1},\ \ \ ~
F_l=a_{l+1,l-1},\ \ \ ~
H_l=a_{l-1,l-1}-a_{l+1,l+1}, 
\end{eqnarray}
where
\begin{eqnarray}
a_{ij}& = & e_{ij}-e_{\jbar\ibar},\ \ \ ~ \nonumber
\ibar=2l+1-i. 
\end{eqnarray}

\subsection{Direct product representation $V^\ep
\otimes V^\ep$\label{subsect:dirprod}}

The direct product representation is constructed in the usual manner 
via the coproduct. For all the algebras $B_l$, $C_l$ and $D_l$, 
$V^\ep~\otimes~V^\ep$ has the following decomposition:
\begin{equation}
V^\ep~\otimes~V^\ep ~~
=~~ V^{2 \ep}~\oplus~V^{\ep+\epsilon_2}~
\oplus~V^0,
\end{equation}
where $V^{\ep+\epsilon_2}$ is the adjoint representation for 
$B_l$ and $D_l$, and $V^{2 \ep}$ is the adjoint representation 
for $C_l$.  In the following we give basis vectors for these 
submodules. The expansion coefficients of these basis
vectors in terms of the $V^\ep$ basis vectors give the quantum
Clebsch-Gordan coefficients (cf. \eqref{exp}).
We have chosen the basis of the vector 
representation to be self-dual, namely $\langle i | j \rangle = 
\delta_{ij}$. Therefore the expansion coefficients of the dual basis
%where dual means dual to the submodule basis, 
are the inverse quantum Clebsch-Gordan coefficients.
We will use a generalized Kronecker delta notation; for 
example, $\delta_{i<j} = 1$ if $i < j$, 0 otherwise.

\subsubsection{$\mathfrak{g}~=~B_l$}
\begin{itemize}
\item Basis for $V^{2 \ep}$:
\begin{eqnarray*}
\omega_{ij} & = & q^{1/2} \mvec{i}{j} + q^{-1/2} \mvec{j}{i}, \ \ \ \ \ \  
i < j \neq \ibar, 
\nonumber\\
\omega_{ii} & = & \mvec{i}{i}, \ \ \ \ \ \  1 \leq i \neq l+1 \leq 2l, 
\nonumber\\
\omega_i & = & q \mvec{i}{\ibar} + q^{-1} \mvec{\ibar}{i}
\nonumber\\ 
&& - ( \mvec{i+1}{\overline{{}i+1}} + 
\mvec{\overline{{}i+1}}{i+1}) \delta_{i < l}
\nonumber\\ 
&& -(q^{1/2} + q^{-1/2}) \mvec{l+1}{l+1} \delta_{il}, \ \ \ \ \ \ i \leq l, 
\end{eqnarray*}
with corresponding dual basis
\begin{eqnarray*}
\omega^i & = & \frac{(q+q^{-1})^{-1}}{D_{l}} \biggl( D_{l-i} 
\sum_{j \leq i} (q^{j} \mvecd{j}{\jbar} + q^{-j} \mvecd{\jbar}{j} ) \\
&& - [i]_q (q^{1/2}+q^{-1/2}) \biggl( \mvecd{l+1}{l+1} \\
&&\ \ \ \ \ \ - \sum_{j > i}^{l} (q^{j-l-1/2} \mvecd{j}{\jbar}
+ q^{l+1/2-j} \mvecd{\jbar}{j}) \biggr) \biggr),
\end{eqnarray*}
where
\begin{eqnarray*}
D_{k} & = & \frac{q^{k+1/2} -q^{-k-1/2}}{q^{1/2}-q^{-1/2}},
\end{eqnarray*}
and the $q$-numbers are given by
\begin{eqnarray*}
[i]_q & = & \frac{q^i-q^{-i}}{q-q^{-1}}.
\end{eqnarray*}

\item Basis for $V^{\epsilon_1 + \epsilon_2}$:
\begin{eqnarray*}
v_{ij} & = & q^{-1/2} \mvec{i}{j} - q^{+1/2} \mvec{j}{i}, \ \ \ \ \ \  
i < j \neq \ibar, \\
v_i & = & \mvec{i}{\ibar} - \mvec{\ibar}{i}\\ 
&& -(q^{-1} \mvec{i+1}{\overline{{}i+1}}
- q \mvec{\overline{{}i+1}}{i+1}) \delta_{i <l} \\
&& + 
(q^{1/2}-q^{-1/2}) \mvec{l+1}{l+1} \delta_{il}, \ \ \ \ \ \ i \leq l,
\end{eqnarray*}
with corresponding dual basis
\begin{eqnarray*}
v^i & = & \frac{(q+q^{-1})^{-1}}{D_l} \biggl( [i]_q (q^{1/2}-q^{-1/2}) 
\biggl( \mvecd{l+1}{l+1} \\
&& + \sum_{j >i}^{l} (q^{j-l-1/2} \mvecd{j}{\jbar} + q^{l+1/2-j} 
\mvecd{\jbar}{j}) \biggr) \\
&& + D_{l-i} \sum_{j \leq i} ( q^{j-1} \mvecd{j}{\jbar} - q^{1-j} 
\mvecd{\jbar}{j}) \biggr),
\end{eqnarray*}
where
\begin{eqnarray*}
D_{k} & = & \frac{q^{k-1/2}+q^{1/2-k}}{q^{1/2}+q^{-1/2}}.
\end{eqnarray*}

\item Basis for $V^0$:
\begin{eqnarray*}
t & = & \sum_{i=1}^{l} \left(q^{i-l-1/2} \mvec{i}{\ibar} + 
q^{-i+l+1/2} \mvec{\ibar}{i}\right) \\
&& + \mvec{l+1}{l+1}.
\end{eqnarray*}
\end{itemize}

\subsubsection{$\mathfrak{g}~=~C_{l}$}

\begin{itemize}
\item Basis for $V^{2 \ep}$:

\begin{eqnarray*}
v_{ij} & = & q^{1/2} \mvec{i}{j} + q^{-1/2} \mvec{j}{i}, \ \ \ \ \ \  
i < j \neq \ibar, \\
v_{ii} & = & \mvec{i}{i}, \ \ \ \ \ \ 1 \leq i \leq 2l,\\
v_i & = & (-1)^i \biggl( q \mvec{i}{\ibar} + q^{-1} \mvec{\ibar}{i} \\
&& + \biggl( \mvec{i+1}{\overline{i+1}} + \mvec{\overline{i+1}}{i+1} 
\biggr) \delta_{i < l} \biggr),\ \ \ \ \ \  1 \leq i \leq l,
\end{eqnarray*}
whose duals are given by
\begin{eqnarray*}
v^i & = & \frac{(q+q^{-1})^{-1}}{D_l} \biggl( D_{l-i} 
\sum_{j \le i} (-1)^j \biggl( q^{j} \mvecd{j}{\jbar} + q^{-j} 
\mvecd{\jbar}{j} \biggr) \\
&& + [i]_q \frac{q-q^{-1}}{q+q^{-1}} \sum_{j > i}^{l} (-1)^j 
\biggl(q^{j-l-1} \mvecd{j}{\jbar}-q^{l+1-j} \mvecd{\jbar}{j} 
\biggr) \biggr),
\end{eqnarray*}
where 
\begin{eqnarray*}
D_{k} & = & \frac{q^{k+1}+q^{-k-1}}{q+q^{-1}}.
\end{eqnarray*}

\item Basis for $V^{\ep+\epsilon_2}$:
\begin{eqnarray*}
\omega_{ij} & = & q^{-1/2} \mvec{i}{j} - q^{1/2} \mvec{j}{i}, \ \ \ \ \ \  
i < j \neq \ibar,\\
\omega_i & = & (-1)^i \biggl(\mvec{i}{\ibar}-\mvec{\ibar}{i} \\
&&  + q^{-1} \mvec{i+1}{\overline{i+1}}-q \mvec{\overline{i+1}}{i+1} 
\biggr), \ \ \ \ \ \  1 \leq i < l,
\end{eqnarray*}
whose corresponding duals are given by
\begin{eqnarray*}
\omega^i & = & \frac{(q+q^{-1})^{-1}}{[l]_q} \biggl( [l-i]_q 
\sum_{j \leq i} (-1)^j \biggl(q^{j-1} \mvecd{j}{\jbar}-q^{1-j} 
\mvecd{\jbar}{j} \biggr) \\
&& + [i]_q \sum_{j > i}^l (-1)^{j-1} \biggl(q^{j-l-1} 
\mvecd{j}{\jbar}-q^{l+1-j} \mvecd{\jbar}{j} \biggr) \biggr).
\end{eqnarray*}

\item Basis for $V^0$:
\begin{eqnarray*}
t & = & \sum_{i=1}^l (-1)^{l-i} \biggl(q^{i-l-1} 
\mvec{i}{\ibar}-q^{l+1-i} \mvec{\ibar}{i} \biggr).
\end{eqnarray*}
\end{itemize}

\subsubsection{$\mathfrak{g}~=~D_l$}
\begin{itemize}
\item Basis for $V^{2 \epsilon_{1}}$:
\begin{eqnarray*}
\omega_{ii} & = & \mvec{i}{i},\ \ \ \ \ \  1 \leq i \leq 2 \mathit l, \\
\mathnormal \omega_{ij} & = & q^{1/2} \mvec{i}{j} + q^{-1/2} 
\mvec{j}{i}, \ \ \ \ \ \  i < j \neq \ibar, \\
\omega_i & = & q \mvec{i}{\ibar} + q^{-1} \mvec{\ibar}{i} \\
&& - \left(\mvec{i+1}{\overline{{}i+1}} + 
\mvec{\overline{{}i+1}}{i+1} \right), \ \ \ \ \ \  1\leq i < l 
\end{eqnarray*}  
with corresponding duals
\begin{eqnarray*}
\omega^{i} & = & \frac{( q+q^{-1} )^{-1}}{[l]_q} \biggl( [l-i]_q 
\sum_{j \leq i} ( q^{j} \mvecd{j}{\jbar} + q^{-j} \mvecd{\jbar}{j} ) \\
&& - [i]_q \sum_{j>i}^{l} (q^{j-l} \mvecd{j}{\jbar} + q^{l-j} 
\mvecd{\jbar}{j} ) \biggr).
\end{eqnarray*}

\item Basis for $V^{\epsilon_1 + \epsilon_2}$: 
\begin{eqnarray*}
v_{ij} & = & q^{-1/2} \mvec{i}{j} - q^{1/2} \mvec{j}{i}, \ \ \ \ \ \  
i < j \neq \ibar, \\
v_i & = & \mvec{i}{\ibar} - \mvec{\ibar}{i} \\
&& - (q^{-1} \mvec{i+1}{\overline{{}i+1}} - 
q \mvec{\overline{{}i+1}}{i+1}) \delta_{i<l}, 
\end{eqnarray*}
with corresponding duals
\begin{eqnarray*}
v^i & = & \frac{(q+q^{-1})^{-1}}{D_l} \biggl( D_{l-i} 
\sum_{j \leq i} (q^{j-1} \mvecd{j}{\jbar} - q^{1-j} 
\mvecd{\jbar}{j}) \\
&& + \frac{[i]_q (q-q^{-1})}{(q+q^{-1})} \sum_{j >i}^{l} 
(q^{j-l} \mvecd{j}{\jbar} + q^{l-j} \mvecd{\jbar}{j}) \biggr),
\end{eqnarray*}
where
\begin{eqnarray*}
D_k & = & \frac{(q^{k-1}+q^{1-k})}{(q+q^{-1})}.
\end{eqnarray*}

\item Basis for $V^0$:
\begin{equation}
t  = \sum_{i=1}^{l} ( q^{i-l} \mvec{i}{\ibar} + q^{l-i} 
\mvec{\ibar}{i}).\nonumber
\end{equation} 
\end{itemize}

Nearly all the elements are in place to construct the quantum 
structure constants. It only remains to identify which vectors 
in the adjoint representation (embedded in $V^\ep
\otimes V^\ep$) correspond to the Cartan subalgebra 
and which to the roots. We have already discussed the weights 
of the basis vectors of the vector representation and therefore we 
can calculate trivially the weights of the basis vectors of the adjoint 
representation. This then fixes the identification between the basis
states of the adjoint representation given above and the basis
states for the quantum Lie algebra chosen in \eqref{qgrad} up to
rescaling. We have choosen the scalings with hindsight so that the
Killing form and the structure constants are as simple as possible.
We set $\xi\, =(q+q^{-1}) \langle t | t \rangle$.  
\newline
For $B_l$ we have
\begin{eqnarray}
X_{\epsilon_i-\epsilon_j} & = & \xi\, v_{i \jbar}, \ \ \ ~ 
X_{-(\epsilon_i-\epsilon_j)} ~ = ~ \xi\, v_{j \ibar}, \nonumber \\
X_{\epsilon_i+\epsilon_j} & = & \xi\, v_{ij}, \ \ \ ~ 
X_{-(\epsilon_i+\epsilon_j)} ~ = ~ \xi\, v_{\jbar \ibar}, \nonumber \\
X_{\epsilon_j} & = & \xi\, v_{j l+1}, \ \ \ ~X_{-\epsilon_j} ~ 
= ~ \xi\, v_{l+1 \jbar}, \nonumber \\
H_j &  = & \xi\, v_j.
\end{eqnarray}
For $C_l$ we have
\begin{eqnarray}
X_{\epsilon_i-\epsilon_j} & = & (-1)^{j-i}\xi\, v_{i \jbar}, \ \ \ ~ 
X_{-(\epsilon_i-\epsilon_j)} ~ = ~ \xi\, v_{j \ibar}, \nonumber \\
X_{\epsilon_i+\epsilon_j} & = & -(-1)^{j-i}\xi\, v_{ij}, \ \ \ ~ 
X_{-(\epsilon_i+\epsilon_j)} ~ = ~ \xi\, v_{\jbar \ibar}, \nonumber \\
X_{2 \epsilon_j} & = & -\xi\, (q+q^{-1})^{1/2} v_{jj}, \ \ \ ~
X_{-2 \epsilon_j} ~ = ~ \xi\, (q+q^{-1})^{1/2} v_{\jbar \jbar}, \nonumber \\
H_i &  = & (-1)^{l+1}\,\xi\, v_i.
\end{eqnarray}
For $D_l$ we have
\begin{eqnarray}
X_{\epsilon_i-\epsilon_j} & = & \xi\, v_{i \jbar}, \ \ \ ~ 
X_{-(\epsilon_i-\epsilon_j)} ~ = ~ \xi\, v_{j \ibar}, \nonumber \\
X_{\epsilon_i+\epsilon_j} & = & \xi\, v_{ij}, \ \ \ ~ 
X_{-(\epsilon_i+\epsilon_j)} ~ = ~ \xi\, v_{\jbar \ibar}, \nonumber \\
H_i &  = &  \xi\, v_i\ \ \ \text{ for }i<l,\\
H_l &=& \xi\, (v_{l-1}+(q+q^{-1})v_l).
\end{eqnarray}
Here $i < j \leq l$, $\ibar = 2l - i +2$ for $B_l$ and 
$\ibar = 2l - i +1$ for $C_l$ and $D_l$.

We could at this point give the structure constants. However some 
of the structure constants are just too complicated. We would 
like to present the structure constants as concisely as possible. 
It turns out that with the introduction of the Killing form on 
our quantum Lie algebras, the unwieldy structure constants 
simplify immensely.

\section{The Killing form\label{sect:kill}}
In this section we introduce an invariant Killing form on every quantum Lie 
algebra. We make the observation that it is an intertwiner for adjoint 
$\otimes$ 
adjoint $\rightarrow$ singlet, where by singlet we mean the trivial
1-dimensional representation $V^0$. This allows us to calculate its
values and
we discover that they are simple $q$-deformation 
of the classical ones. To begin, we define the quantum 
analogue to the classical Killing form as follows:

\begin{definition} The \dem{quantum Killing form} is the map 
$\mathfrak{B}: \mathfrak{L}_h(\mathfrak{g}) \hat{\otimes} 
\mathfrak{L}_h(\mathfrak{g}) \rightarrow \mathbb C[[h]]$  
given by
\begin{equation} \label{kill}
\mathfrak{B}(a,b) = Tr_\Psi  (a\,b\,u).
\end{equation}
\end{definition}
Here $u$ is the the element of $U_h(\mathfrak{g})$ satisfying the 
properties $u\,a\,u^{-1} = S^2(a) \quad\forall  a \in U_h(\mathfrak{g})$ 
and $\Delta(u)= u \otimes u$. The $Tr_{\Psi}$ denotes the trace 
over the adjoint representation. This definition for the quantum 
Killing form reduces to that of the classical Killing form in the classical 
limit. It obviously exists and is non-degenerate because degeneracy 
would spoil the non-degeneracy of the classical Killing form. 
The definition is motivated by the following:
\begin{proposition} 
\label{prop:adinv}
The Killing form is ad-invariant, i.e.
\begin{eqnarray}
\label{adinv}
\mathfrak{B}([a , b]_h,c) & = & \mathfrak{B}(a,[b , c]_h), \ \ \ \ \ \  
\forall a,b,c \in \mathfrak{g}_h
\end{eqnarray}
\end{proposition}
We will prove this at the end of this section.

The 
following proposition, whose proof is trivial, informs us that unlike 
the classical Killing form, the quantum Killing form is not symmetric.
\begin{proposition} The quantum Killing form $\mathfrak{B}$ is a 
non-degenerate, bilinear, nonsymmetric form. Symmetry is replaced by 
the relation
\begin{equation}
\label{nsymm}
\mathfrak{B}(a,b) = \mathfrak{B}(b,S^2(a)). 
\end{equation}
\end{proposition}
The square of the antipode $S^2$ acts on the basis elements of 
$\mathfrak{L}(\mathfrak{g})$ by multiplication by a power of $q$. 
Therefore $S^2(\mathfrak{L}(\mathfrak{g})) 
\subset \mathfrak{L}(\mathfrak{g})$ and \eqref{nsymm} makes sense. 

The calculation of the Killing form is made simple when we realize 
that the Killing 
form is an intertwiner for adjoint 
$\otimes$ adjoint $\rightarrow$ singlet.

\begin{proposition} The Killing form defined in \eqref{kill} is 
an intertwiner from adjoint $\otimes$ adjoint 
$\rightarrow$ singlet, i.e.,
\begin{equation}
\label{kprop}
\mathfrak{B}([c_{(1)} , a]_h,[c_{(2)} , b]_h)
=\epsilon(c) \mathfrak{B}(a,b), \ \ \ ~ 
\forall a,b,c \in \mathfrak{L}_h(\mathfrak{g})
\end{equation}
\end{proposition}
\begin{proof} From the definition of the Killing form and the Lie bracket, 
the left hand side becomes
\begin{eqnarray*}
Tr_{\Psi}(c_{(1)}a S(c_{(2)}) c_{(3)} b S(c_{(4)}) u) & 
= & Tr_{\Psi}(c_{(1)} ab S(c_{(2)})u) \\
& = & Tr_{\Psi}(c_{(1)} ab u S^{-1}(c_{(2)})) \\
& = & Tr_{\Psi}(S^{-1}(c_{(2)}) c_{(1)} ab u) \ \ \ ~~
\mbox{by cyclicity} \\
& = & Tr_{\Psi}(\epsilon(c) ab u) 
\end{eqnarray*}
\end{proof}

Such an intertwiner is unique up to rescaling. To
see this let
$h,g : V^{\Psi} \otimes V^{\Psi} \rightarrow V^0$
be any two intertwiners.
We represent them by thin resp. thick lines as follows:
\begin{equation}
% Drawing generated by LaTeX-CAD 1.7 - requires latexcad.sty 
% (c) 1996 John Leis leis@usq.edu.au 
\begin{picture}(108,26)
\thinlines
\drawarc{18.0}{16.0}{12.0}{6.28}{3.14}
\drawrighttext{8.0}{14.0}{$h~=$}
\drawrighttext{48.0}{14.0}{$h^{-1}~=$}
\Thicklines
\drawarc{102.0}{16.0}{12.0}{6.28}{3.14}
\thinlines
\drawrighttext{92.0}{14.0}{$g~=$}
\drawarc{58.0}{10.0}{12.0}{3.14}{6.28}
\end{picture}
\end{equation}
We can compose either $h$ or $g$ with $h^{-1}$ to obtain
two intertwiners from $V^\Psi$ to itself. 
\begin{equation}
% Drawing generated by LaTeX-CAD 1.7 - requires latexcad.sty 
% (c) 1996 John Leis leis@usq.edu.au 
\begin{picture}(154,22)
\thinlines
\drawarc{10.0}{10.0}{12.0}{6.284}{3.142}
\drawarc{22.0}{10.0}{12.0}{3.14}{6.28}
\drawpath{4.0}{10.0}{4.0}{16.0}
\drawpath{28.0}{10.0}{28.0}{4.0}
\drawpath{48.0}{16.0}{48.0}{4.0}
\Thicklines
\drawarc{74.0}{10.0}{12.0}{6.28}{3.14}
\thinlines
\drawarc{86.0}{10.0}{12.0}{3.14}{6.28}
\drawpath{68.0}{10.0}{68.0}{16.0}
\drawpath{92.0}{10.0}{92.0}{4.0}
\drawcenteredtext{38.0}{10.0}{$\propto$}
\drawcenteredtext{58.0}{10.0}{$\propto$}
\drawarc{120.0}{12.0}{12.0}{6.28}{3.14}
\Thicklines
\drawarc{148.0}{12.0}{12.0}{6.28}{3.14}
\thinlines
\drawcenteredtext{134.0}{10.0}{$\propto$}
\drawcenteredtext{104.0}{10.0}{$\Longrightarrow$}
\end{picture}
\end{equation}
By Schur's lemma they both have to be proportional to the identity.
Thus $h$ and $g$ have to be proportional.

One can also define a form on all of $\uqg$ by equation \eqref{kill}.
This gives the Rosso form \cite{Ros90} which is the unique form
on $\uqg$ for which \eqref{kprop} is valid for any $a,b,c\in\uqg$.
 
We calculate an intertwiner $B$ from adjoint $\otimes$ adjoint 
$\rightarrow$ singlet using the
Clebsch-Gordan coefficients from section \ref{subsect:dirprod}
according to the formula
\begin{equation}\label{killeq}
% Drawing generated by LaTeX-CAD 1.7 - requires latexcad.sty 
% (c) 1996 John Leis leis@usq.edu.au 
\begin{picture}(110,34)
\thinlines
\drawarc{34.0}{18.0}{16.0}{6.28}{3.14}
\drawrighttext{14.0}{14.0}{$B(a,b)~=$}
\drawcenteredtext{58.0}{14.0}{$=~\frac{1}{N}$}
\drawrighttext{24.0}{18.0}{$\Psi$}
\drawcenteredtext{26.0}{22.0}{$a$}
\drawcenteredtext{42.0}{22.0}{$b$}
\drawlefttext{44.0}{18.0}{$\Psi$}
\drawpath{78.0}{24.0}{78.0}{20.0}
\drawarc{86.0}{16.0}{24.0}{6.28}{3.14}
\drawarc{78.0}{16.0}{8.0}{3.14}{6.28}
\drawarc{86.0}{16.0}{8.0}{6.28}{3.14}
\drawarc{94.0}{16.0}{8.0}{3.14}{6.28}
\drawpath{94.0}{24.0}{94.0}{20.0}
\drawlefttext{100.0}{16.0}{$\ep$}
\drawrighttext{72.0}{16.0}{$\ep$}
\drawlefttext{96.0}{22.0}{$\Psi$}
\drawrighttext{76.0}{22.0}{$\Psi$}
\drawcenteredtext{78.0}{28.0}{$a$}
\drawcenteredtext{94.0}{28.0}{$b$}
\drawcenteredtext{86.0}{16.0}{$\ep$}
\end{picture}
\end{equation}
where $N = (q+q^{-1})^{3}$ for $B_l$ and $D_l$, and 
$N=-(q+q^{-1})^{3}$ for $C_l$. By the above propositions 
$B$ is proportional to the quantum Killing form $\mathfrak{B}$
defined in \ref{kill}. 
%The relationship between them is summarized below:
%\begin{eqnarray}
%\mathfrak{B}(a,b) & = & \bigl(\frac{2(q^{l+3/2}+q^{-l-3/2})
%(q+q^{-1})}{q^{l-1/2}+q^{-l+1/2}}-(q^3+q^{-3}) \nonumber  \\
%& + & \frac{(q^{2l-2}-q^{-2l+2})(q+q^{-1})}{q-q^{-1}}\bigr) 
%B(a,b) \ \ \ ~\mbox{   for $B_{l}$,} \nonumber \\
%\mathfrak{B}(a,b) & = & \bigl( \bigl(\frac{q^{3l+3}+q^{-3l-3}+
%2(q^{l+3}+q^{-l-3})}{q^{l+1}+q^{-l-1}}-1 \bigr)(q+q^{-1}) \nonumber \\
%& + & (q^{2l-2}+q^{-2l+2})(q^2+q^{-2})(q+q^{-1})   \\
%& + & (\frac{q^{2l-1}-q^{-2l+1}}{q^3-q^{-3}}-1)(q^3+q^{-3})\bigr) 
%B(a,b)\ \ \ ~\mbox{   for $C_{l}$,} \nonumber \\
%\mathfrak{B}(a,b) & = & \bigl( \frac{q^{4-l}+q^{l-4}}
%{q^{l-1}+q^{1-l}} + \frac{q^{2l-3}-q^{3-2l}-q+q^{-1}}
%{q^3-q^{-3}}(q^3+q^{-3})\bigr) B(a,b) \ \ \ ~
%\mbox{   for $D_{l}$} \nonumber. 
%\end{eqnarray}

The Killing form on the roots has the simple form 
\begin{eqnarray}
\label{kroot}
B(X_{\alpha},X_{\beta}) & = & q^{-\rho \cdot \alpha}~ 
\delta_{\alpha,-\beta}, \ \ \ \ \ \  \forall \alpha \in R. 
\end{eqnarray}
The Killing form on the Cartan subalgebra  for 
the algebras $A_l$, $B_l$, $C_l$ and $D_l$ respectively, 
is given by 
 
\begin{equation} \label{kill0}
(A_l)_h:\ \ \ ~~B(H_i,H_j) ~=~
\begin{pmatrix} [2]_q & -1 & & & & & & & \\
                -1 & [2]_q & -1 & & & & & & \\
                 & & & & & & & & \\
                 & & & & & & & & \\
                 & & & & & & & & \\
                 & & & & & & -1 & [2]_q & -1 \\
                 & & & & & & & -1 & [2]_q
\end{pmatrix}
\end{equation}
\begin{equation} \label{kill1}
(B_l)_h:\ \ \ ~~B(H_i,H_j) ~=~
\begin{pmatrix} [2]_q & -1 & & & & & & & \\
                -1 & [2]_q & -1 & & & & & & \\
                 & & & & & & & & \\
                 & & & & & & & & \\
                 & & & & & & & & \\
                 & & & & & & -1 & [2]_q & -1 \\
                 & & & & & & & -1 & 1
\end{pmatrix}
\end{equation}
\begin{equation}
\label{kill2}
(C_l)_h:\ \ \ ~~ B(H_i,H_j) ~=~
\begin{pmatrix} [2]_q & -1 & & & & & & & \\
                -1 & [2]_q & -1 & & & & & & \\
                 & & & & & & & & \\
                 & & & & & & & & \\
                 & & & & & & & & \\
                 & & & & & & -1 & [2]_q & -1 \\
                 & & & & & & & -1 & \frac{q^2+q^{-2}}{q+q^{-1}}
\end{pmatrix}
\end{equation}
\begin{equation}
\label{kill3}
(D_l)_h:\ \ \ ~~ B(H_i,H_j) ~=~
\begin{pmatrix} [2]_q & -1 & & & & & & & \\
                -1 & [2]_q & -1 & & & & & & \\
                 & & & & & & & & \\
                 & & & & & & & & \\
                 & & & & & & & & \\
                 & & & & & & -1 & [2]_q & -1 & -1 \\
                 & & & & & & & -1 & [2]_q & 0 \\
                 & & & & & & & -1 & 0 & [2]_q
\end{pmatrix}
\end{equation}
These expressions are surprisingly simple $q$-deformations of the
classical ones. In particular $B(H_i,H_j)\neq 0$ if and only if 
$\alpha_i\cdot\alpha_j\neq 0$. The Killing form for $(A_l)_h$
was determined in \cite{Del95}.

We introduce the notation $B(H_i,H_j)=B_{ij}, 
B(X_{\alpha},X_{-\alpha})=B_{\alpha,-\alpha}$. We define 
$B^{ij}$ such that $B_{ij} B^{jk}= \delta_{i}^k$ 
and  $B^{-\alpha,\alpha}$ such that $B_{\alpha,-\alpha}
B^{-\alpha,\alpha}=1$. Or, equivalently, introducing composite
indices $p,q,r=\{i,\alpha\}$,  
\begin{equation}
\setlength{\unitlength}{0.65mm}
% Drawing generated by LaTeX-CAD 1.7 - requires latexcad.sty 
% (c) 1996 John Leis leis@usq.edu.au 
\begin{picture}(186,32)
\thinlines
\drawrighttext{56.0}{16.0}{$B^{pq}\,=$}
\drawarc{26.0}{18.0}{12.0}{6.28}{3.14}
\drawcenteredtext{88.0}{16.0}{s.t.}
\drawcenteredtext{20.0}{22.0}{$p$}
\drawcenteredtext{32.0}{22.0}{$q$}
\drawrighttext{14.0}{16.0}{$B_{pq}\,=$}
\drawpath{126.0}{22.0}{126.0}{16.0}
\drawarc{66.0}{16.0}{12.0}{3.14}{6.28}
\drawcenteredtext{60.0}{12.0}{$p$}
\drawcenteredtext{72.0}{12.0}{$q$}
\drawpath{150.0}{16.0}{150.0}{10.0}
\drawarc{132.0}{16.0}{12.0}{6.28}{3.14}
\drawarc{144.0}{16.0}{12.0}{3.14}{6.28}
\drawpath{164.0}{22.0}{164.0}{10.0}
\drawcenteredtext{158.0}{16.0}{$=$}
\drawlefttext{168.0}{16.0}{$=~\delta^p_r$}
\drawcenteredtext{126.0}{26.0}{$p$}
\drawcenteredtext{164.0}{26.0}{$p$}
\drawrighttext{122.0}{16.0}{$B_{pr}B^{rq}\,=$}
\drawcenteredtext{150.0}{6.0}{$q$}
\drawcenteredtext{164.0}{6.0}{$q$}
\end{picture}
\setlength{\unitlength}{0.8mm}
\end{equation}
We give the formulas for $B^{ij}$ with $i\leq j$, the remaining
ones are obtained by symmetry, $B^{ji}=B^{ij}$.
\begin{align}\label{killinvf}
(A_l)_h:&&B^{ij}=\,
&[i]_q\,\frac{q^{l-j+1}-q^{j-l-1}}{q^{l+1}-q^{-l-1}},
\\
(B_l)_h:&&B^{ij}=\,
&[i]_q\,\frac{q^{l-j-1/2}+q^{j-l+1/2}}{q^{l-1/2}+q^{-l+1/2}},
\\
(C_l)_h:&&B^{ij}=\,
&[i]_q\,\frac{q^{l-j+1}+q^{j-l-1}}{q^{l+1}+q^{-l-1}},
\\
(D_l)_h:&&B^{ij}=\,
&[i]_q\,\frac{q^{l-j-1}+q^{j-l+1}}{q^{l-1}+q^{-l+1}},\ \ \ 
i<l-1,j<l,
\nonumber\\
&&B^{l-1,l-1}=B^{l,l}=\,
&[l]_q\frac{(q+q^{-1})^{-1}}{q^{l-1}+q^{-l+1}},
\nonumber\\
&&B^{l-1,l}=\,
&[l-2]_q\frac{(q+q^{-1})^{-1}}{q^{l-1}+q^{-l+1}},
\nonumber\\
&&B^{i,l}=&B^{i,l-1}.\label{killinvl}
\end{align}

We now give the proof of Proposition \ref{prop:adinv}, i.e., we
will show that
\begin{equation}
\label{adinv2}
B([a , b]_h,c)  =  B(a,[b , c]_h), \ \ \ \ \ \  
\forall a,b,c \in \mathfrak{g}_h
\end{equation}
\begin{proof}
Because the intertwiner from adjoint to adjoint $\otimes$ adjoint
is unique up to rescaling we have 
\begin{equation}\label{gr}
% Drawing generated by LaTeX-CAD 1.7 - requires latexcad.sty 
% (c) 1996 John Leis leis@usq.edu.au 
\begin{picture}(76,36)
\thinlines
\drawarc{12.0}{24.0}{12.0}{6.28}{3.14}
\drawpath{12.0}{18.0}{12.0}{10.0}
\drawcenteredtext{12.0}{6.0}{$r$}
\drawcenteredtext{6.0}{30.0}{$p$}
\Thicklines
\drawcenteredtext{18.0}{30.0}{$q$}
\thinlines
\drawpath{6.0}{24.0}{6.0}{26.0}
\drawpath{18.0}{24.0}{18.0}{26.0}
\Thicklines
\drawcenteredtext{32.0}{18.0}{$=~~A~~$}
\thinlines
\drawarc{56.0}{20.0}{8.0}{6.28}{3.14}
\drawarc{50.0}{16.0}{12.0}{6.28}{3.14}
\drawarc{64.0}{20.0}{8.0}{3.14}{6.28}
\drawpath{48.0}{18.0}{48.0}{18.0}
\drawpath{48.0}{18.0}{48.0}{18.0}
\drawpath{52.0}{20.0}{52.0}{26.0}
\drawpath{44.0}{16.0}{44.0}{26.0}
\drawpath{68.0}{20.0}{68.0}{10.0}
\drawcenteredtext{44.0}{30.0}{$p$}
\Thicklines
\drawcenteredtext{52.0}{30.0}{$q$}
\thinlines
\drawcenteredtext{68.0}{6.0}{$r$}
\end{picture}
\end{equation}
The constant of proportionality 
$A\in \mathbb C[[h]]$
can be determined by setting $p=r=X_{\alpha}$ and  $q=H_i$. Doing this 
we obtain the relation
\begin{eqnarray*}
-r_{\alpha} ~=~A~ B_{\alpha,-\alpha}~ l_{-\alpha}(H_i) ~
B^{-\alpha,\alpha}.
\end{eqnarray*}
Using $B_{\alpha,-\alpha}
B^{-\alpha,\alpha}=1$ and the relation
$l_{-\alpha} = -r_{\alpha}$ from \eqref{qant1} and \eqref{qant2}
we see that $A=1$. 

We now write the left hand side
of \eqref{adinv2} in graphical form and use the above identity 
\eqref{gr} to 
manipulate the diagram to the right hand side of \eqref{adinv2}. 
\begin{equation}
% Drawing generated by LaTeX-CAD 1.7 - requires latexcad.sty 
% (c) 1996 John Leis leis@usq.edu.au 
\begin{picture}(100,24)
\thinlines
\drawarc{6.0}{16.0}{8.0}{6.28}{3.14}
\drawarc{12.0}{12.0}{12.0}{6.28}{3.14}
\drawpath{18.0}{12.0}{18.0}{16.0}
\drawcenteredtext{26.0}{12.0}{$=$}
\drawcenteredtext{74.0}{12.0}{$=$}
\drawarc{88.0}{12.0}{12.0}{6.28}{3.14}
\drawarc{94.0}{16.0}{8.0}{6.28}{3.14}
\drawpath{82.0}{12.0}{82.0}{16.0}
\drawpath{40.0}{12.0}{40.0}{12.0}
\drawpath{40.0}{12.0}{40.0}{12.0}
\drawarc{46.0}{14.0}{8.0}{6.28}{3.14}
\drawarc{40.0}{10.0}{12.0}{6.28}{3.14}
\drawarc{54.0}{14.0}{8.0}{3.14}{6.28}
\drawpath{38.0}{12.0}{38.0}{12.0}
\drawpath{38.0}{12.0}{38.0}{12.0}
\drawpath{42.0}{14.0}{42.0}{18.0}
\drawpath{34.0}{10.0}{34.0}{18.0}
\drawarc{62.0}{14.0}{8.0}{6.28}{3.14}
\drawpath{66.0}{14.0}{66.0}{18.0}
\end{picture}
\end{equation}
\end{proof}

Besides the Killing form there exists another natural bilinear form
$\rp{}{}$ on the quantum Lie algebras. This is the unique (up to
rescaling) non-degenerate bilinear form which satisfies
\begin{equation}
\rp{[a,b]}{c}=\rp{b}{[a^\dagger,c]}\ \ \ ~\forall a,b,c\in
\mathfrak{L}_h(\mathfrak{g}),
\end{equation}
where $\dagger$ denotes the algebra antiautomorphism of
$U_h(\lie)$ defined by
\begin{equation}
(x_i^\pm)^\dagger=x^\mp_i,\ \ \ \ \ \ ~h_i^\dagger=h_i.
\end{equation}
(If we define $\dagger$ to be anti-linear then the form $\rp{}{}$ has 
to be sesquilinear.)
On the quantum Lie algebra generators it acts as
\begin{equation}
X_\alpha^\dagger=q^{-\rho\cdot\alpha}\,X_{-\alpha},\ \ \ \ \ \ ~
H_i^\dagger=H_i.
\end{equation}
The form $\rp{}{}$ is symmetric and is given by
\begin{equation}
\rp{X_\alpha}{X_\beta}=\delta_{\alpha\beta},\ \ \ ~~
\rp{H_i}{H_j}=B_{ij},\ \ \ ~~\rp{X_\alpha}{H_i}=\rp{H_i}{X_\alpha}=0.
\end{equation}

\section{The quantum structure constants\label{sect:results}}
We only give the minimal set of the structure constants. The remaining 
ones can be calculated using the symmetry properties, see equations
(\ref{qant1})-(\ref{qant5}). 
Thus we give only the left quantum roots $l_\alpha$ for positive 
$\alpha$. From these those for negative $\alpha$ as well as the
$r_\alpha$ and $g_\alpha{}^k$ can be
obtained using eqs. \eqref{qant1}, \eqref{qant2} and \eqref{qant5}.
The structure constants $f_{ij}{}^k$ for the Cartan subalgebra simplify 
dramatically if one lowers the last index using the Killing form, i.e.,
$f_{ijk} := f_{ij}{}^{m} B(H_m,H_k)$. Since the Killing 
form is symmetric on the Cartan subalgebra it is easy to see, using 
\eqref{qant1},\eqref{qant2} and 
\eqref{qant4}, that $f_{ijk}$ is completely symmetric. We also 
observed that $f_{ijk} \neq 0$
iff $\alpha_i\cdot\alpha_j\neq 0$. This implies in particular that
the structure constants $f_{ijk}$ are non-zero only if at least two 
of the indices are the same (This is so because $f_{ijk}$ is non-zero only 
if $f_{ijk}$, $f_{ikj}$ and $f_{jki}$ are non-zero.
In other words, $f_{ijk}$ is non-zero if there exist three simple roots such 
that $\alpha_i\cdot\alpha_j$, $\alpha_i\cdot\alpha_k$
and $\alpha_j\cdot\alpha_k$ are non-zero. There exist no three distinct
simple roots that satisfy this requirement.)
So below we give only a few $f$'s, all others can be obtained by
symmetry or are zero.
We give only enough of the $N_{\alpha,\beta}$ so that the others can be
obtained using \eqref{qant1} -- \eqref{qant3}.

\subsection{The quantum Lie algebra ($B_l$)$_h$}

The quantum roots are
\begin{eqnarray}\label{ba}
l_{\epsilon_j-\epsilon_k} (H_i) & = &  
q^{l-i-1/2} \delta_{ij} - q^{i-l+1/2} \delta_{ik} 
- q^{l-i-5/2} \delta_{i+1,j} + q^{i-l+5/2} \delta_{i+1,k}, \nonumber\\
l_{\epsilon_j+\epsilon_k} (H_i) & = &   
q^{l-i-1/2} \delta_{ij} + q^{l-i+3/2} \delta_{ik} 
- q^{l-i-5/2} \delta_{i+1,j} - q^{l-i-1/2} \delta_{i+1,k}, \nonumber\\
l_{\epsilon_k} (H_i) & = &    
q^{l-i-1/2} \delta_{ik}-q^{l-i-5/2} \delta_{i+1,k} 
+ (q^{3/2}-q^{1/2}) \delta_{il} ,
\end{eqnarray}
where $j < k \leq l$, $i \leq l$.
The structure constants for the Cartan subalgebra are
\begin{eqnarray}\label{bf}
f_{iii} & = & (q^{l-i-3/2}+q^{i-l+3/2})(q^2-q^{-2}), \nonumber\\
f_{lll} & = & (q+q^{-1})(q^{1/2}-q^{-1/2})-(q^{1/2}-q^{-1/2}), \nonumber\\
f_{i\pm 1,i\pm 1,i} & = & \mp\ (q^{l-i-3/2}-q^{i-l+3/2}), 
\end{eqnarray}
The remaining structure constants are determined by
\begin{eqnarray}\label{bn}
N_{\epsilon_i-\epsilon_j,\epsilon_k-\epsilon_m} & = &  q^{l-k-1} 
\delta_{jk} ~-~ q^{m-l+1} \delta_{im}, \nonumber\\
N_{\epsilon_i-\epsilon_j,\epsilon_k+\epsilon_m} & = &  q^{l-k-1} 
\delta_{jk} ~+~  (q^{l-m+1}  \delta_{i > k} ~-~ q^{l-m} \delta_{i < k}) 
\delta_{jm}, \nonumber\\
N_{\epsilon_i-\epsilon_j,\epsilon_m} & = &  q^{l-m-1} \delta_{jm}, \nonumber\\
N_{\epsilon_j,\epsilon_m} & = &  q^{1/2} \delta_{j > m}~-~ q^{-1/2} 
\delta_{j < m},
\end{eqnarray}
where $i< j \leq l$ and $k < m \leq l$.

\subsection{The quantum Lie algebra (C$_l$)$_h$}

The quantum roots are
\begin{eqnarray}\label{ca}
l_{\epsilon_j-\epsilon_k}(H_i) & = &
q^{l-i+3} \delta_{ij}- q^{i-l-3}\delta_{ik} 
- q^{l-i+1}\delta_{i+1,j} + q^{i-l-1} \delta_{i+1,k}, \nonumber\\
l_{\epsilon_j+\epsilon_k}(H_i) & = & 
q^{l-i+3} \delta_{ij} + q^{l-i+1} \delta_{ik} 
- q^{l-i+1} \delta_{i+1,j} - q^{l-i-1} \delta_{i+1,k},\nonumber\\
l_{2 \epsilon_k} (H_i) & = & 
(q+q^{-1}) (q^{l-i+2}\delta_{ik} - q^{l-i} \delta_{i+1,k}), 
\end{eqnarray}
where $j < k \leq l$, $i \leq l$.  
The structure constants for the Cartan subalgebra are
\begin{eqnarray}\label{cf}
f_{iii} & = & (q^{2}-q^{-2})(q^{l-i+2}+q^{i-l-2}), \nonumber\\
f_{i\pm 1,i\pm 1,i}  & = & \mp\  (q^{l-i+2}-q^{i-l-2}).
\end{eqnarray}
The remaining structure constants are determined by
\begin{eqnarray}\label{cn}
N_{\epsilon_i-\epsilon_j,\epsilon_k-\epsilon_m} & = &  (-1)^{l-i}
q^{i-l-5/2} \delta_{im} ~-~ (-1)^{l-j} q^{l-j+5/2} \delta_{jk}, \nonumber\\ 
N_{\epsilon_i-\epsilon_j,\epsilon_k+\epsilon_m} & = & - (-1)^{l-j} \bigl\{
q^{l-j+5/2} \delta_{jk} + \left(q^{l-j+3/2} \delta_{i < k}\right.\nonumber\\
&&\ \ \ ~\left.+(-1)^{j-i}\,(q+q^{-1})^{1/2}  q^{l-j+1} \delta_{ik}~+~  
q^{l-j+1/2} \delta_{i > k}\right) \delta_{jm} \bigr\},\nonumber\\ 
N_{\epsilon_i-\epsilon_j,2 \epsilon_m} & = & - (-1)^{l-j} (q+q^{-1})^{1/2}
q^{l-j+2} \delta_{jm}, \nonumber\\
N_{\epsilon_i-\epsilon_j,-2 \epsilon_m} & = &  (-1)^{l-j} (q+q^{-1})^{1/2}
q^{i-l-1} \delta_{im},
\end{eqnarray}
where $i < j \leq l$ and $k < m \leq l$.

\subsection{The quantum Lie algebra ($D_l$)$_h$}

The quantum roots are
\begin{eqnarray}\label{da}
l_{\epsilon_j-\epsilon_k} (H_i) & = & (q^{l-i-1} \delta_{ij} - 
q^{i-l+1} \delta_{ik} - q^{l-i-3} \delta_{i+1,j} + 
q^{i-l+3} \delta_{i+1,k} ) \delta_{i < l} \nonumber\\
&& + (- \delta_{ik}+\delta_{i-1,j}-\delta_{i-1,k})\delta_{il}, \nonumber\\
l_{\epsilon_{j} + \epsilon_{k}} (H_i) & = &   
(q^{l-i-1} \delta_{ij} + q^{l-i+1} \delta_{ik} - 
q^{l-i-3} \delta_{i+1,j} - q^{l-i-1} \delta_{i+1,k} ) \delta_{i < l} \nonumber\\
&& + (q^2 \delta_{ik}+\delta_{i-1,j}+q^2 \delta_{i-1,k})\delta_{il},  
\end{eqnarray}
where $j < k \leq l$, $i \leq l$.
The structure constants for the Cartan subalgebra are
\begin{eqnarray}\label{df}
f_{iii} & = & (q^{2}-q^{-2})(q^{l-i-2}+q^{i-l+2}), \nonumber\\
f_{j+1,j+1,j} & = & - (q^{l-j-2}-q^{j-l+2}), \nonumber\\
f_{i-1,i-1,i} & = & + (q^{l-i-2}-q^{i-l+2}), 
\end{eqnarray} 
where $i<l, j<l-1$. Because of the Dynkin diagram automorphism 
$\tau$ which interchanges
$H_l$ and $H_{l-1}$ we don't need to give the structure
constants involving $H_l$, they are equal to those
involving $H_{l-1}$. The $f$'s involving both $l$ and $l-1$ are
zero. The remaining structure constants are determined by
\begin{eqnarray}\label{dn}
N_{\epsilon_i-\epsilon_j,\epsilon_k-\epsilon_m} & = &  
q^{l-j-3/2} \delta_{jk} ~-~ q^{i-l+3/2} \delta_{im}, \nonumber\\
N_{\epsilon_i-\epsilon_j,\epsilon_k+\epsilon_m} & = &  
q^{l-j-3/2} \delta_{jk} ~-~ (q^{l-j-1/2} \delta_{i < k} ~-~ 
q^{l-j+1/2} \delta_{i > k}) \delta_{jm},
\end{eqnarray} 
where $i<j \leq l$ and $k<m \leq l$.

\subsection{The quantum Lie algebras ($A_l$)$_h$}

For completeness we also give the structure constants for the quantum 
Lie algebras associated to $\lie=A_l$ which were determined by a 
different method in \cite{Del95}, see also \cite{Del97}. There is a
family of quantum Lie algebras $(A_l)_h(\chi)$ depending on a
parameter $\chi$. This is due to the fact that in the case of $A_l$
the adjoint representation appears in adjoint $\otimes$ adjoint with 
multiplicity two. The parameter $\chi$ can be written as a fraction
$\chi=s/t$ with $s,t\in\ch$ and with the retriction that $(s+t)^{-1}
\in\ch$. $A_l$ is isomorphic to $sl_{l+1}$.

The quantum roots are
\begin{align}\label{aa}
l_{\epsilon_j-\epsilon_k}(H_i)=&(q^{1-i}\delta_{ij}-q^{-1-i}\delta_{i+1,j})
(s+t\,q^{l+1})\nonumber\\
&-(q^{i-1}\delta_{ik}-q^{i+1}\delta_{i+1,k})(s+t\,q^{-l-1}),
\end{align}
where $j \neq k \leq l+1$, $i \leq l$.
The structure constants for the Cartan subalgebra are
\begin{align}\label{af}
f_{iii} = & s (q^{2}-q^{-2})(q^{-i}+q^{i}) 
+t (q^{2}-q^{-2})(q^{l-i+1}+q^{i-l-1}), \nonumber\\
f_{i\pm 1,i\pm 1,i}  = & \mp s (q^{-i}-q^{i})
\mp t (q^{l-i+1}-q^{i-l-1}), 
\end{align}
where $i\leq l$. And finally
\begin{equation}\label{an}
N_{\epsilon_i-\epsilon_j,\epsilon_k-\epsilon_m}=
q^{1/2-j}(s+t\,q^{l+1})\delta_{jk}-
q^{i-1/2}(s+t\,q^{-l-1})\delta_{im},
\end{equation}
where $i\neq j \leq l+1$ and $k\neq m \leq l+1$.

\subsection{The structure of the Cartan subalgebra}

One of the novel features of quantum Lie algebras is the non-vanishing
of the quantum Lie bracket between elements of the Cartan subalgebra.
Our explicit results for the corresponding structure constants
$f_{ijk}$ given above have lead us to the following observation:
\begin{equation}\label{fr}
f_{iij}=B_{ij}\,\left(l_{\alpha_j}(H_i)-
r_{\alpha_j}(H_i)\right).
\end{equation}
In other words, the Lie brackets of the Cartan subalgebra elements are
given by the amount of the split between left and right quantum roots.
In the following we will abbreviate this $q$-antisymmetric combination
of left and right roots by $a_\alpha:=l_\alpha-r_\alpha$. The Lie bracket
relations are then
\begin{equation}
[H_i,H_j]_h=B_{ij}\,\left(a_{\alpha_i}(H_j)\,H^j+
a_{\alpha_j}(H_i)\,H^i\right)
\end{equation}

\subsection{The quantum root space}

The quantum root space $\cart^*$ is the
dual space to the Cartan subalgebra, i.e., it is the space of linear
functionals on $\cart$ with values in $\ch$. The left and right
quantum roots $l_\alpha$ and $r_\alpha$ are particular
elements of $\cart^*$.

The Killing form $B$ on $\cart$ provides a natural pairing between
elements $H\in\cart$ and linear functionals $v_H\in\cart^*$ defined
by
\begin{equation}
v_H(H'):=B(H,H')\ \ \ ~\forall H'\in\cart.
\end{equation}
Let $v_i:=v_{H_i}$ be the elements of $\cart^*$ dually paired with
the generators $H_i$ of $\cart$. 
Then $\{v_i\}_{i=1,\dots,\text{rank}\lie}$ is
a basis for $\cart^*$. With our choice of the $H_i$ the 
$v_i$ are proportional to
the $q$-symmetric combination of the simple left and right
quantum roots,
\begin{equation}\label{br}
v_i(H_j):=B(H_i,H_j)=\frac{1}{\xi_i}\,\left(l_{\alpha_i}(H_j)+
r_{\alpha_i}(H_j)\right).
\end{equation}
The factors of proportionality $\xi_i$ are given by
\begin{align}
B_l:\ \ \ ~~&\xi_i=q^{i-l+3/2}+q^{l-i-3/2},\\
C_l:\ \ \ ~~&\xi_i=q^{i-l-2}+q^{l-i+2}~~\text{for }i<l,\ \ \ ~~
\xi_l=(q+q^{-1})^2,\\
D_l:\ \ \ ~~&\xi_i=q^{i-l+2}+q^{l-i-2}~~\text{for }i<l,\ \ \ ~~
\xi_l=\xi_{l-1}=(q+q^{-1}).
\end{align}
These were chosen so as to make the structure constants as simple
as possible.

There is a natural inner product on $\cart^*$ given by
\begin{equation}
\rp{v_H}{v_{H'}}:=B(H,H').
\end{equation}
On our basis this gives $\rp{v_i}{v_j}=B_{ij}$. Thus we have
chosen our basis vectors to all have length squared equal to
$[2]_q=(q+q^{-1})$ except for $v_l$ for $B_l$ and $C_l$. The
simple quantum roots $l_{\alpha_i}$ on the other hand all have
different lenghts, quite unlike the classical simple roots.

\section{Some more relations between the structure constants
\label{sect:rels}}
We have already observed that the product of any number of intertwiners 
is still an intertwiner. This observation is indeed a very powerful one. 
Below we will use it to derive some interesting relations between the 
structure constants. 

The calculations leading to our results for the structure constants
presented in the previous section were rather lengthy. It is
therefore very important to have powerful checks on the correctness
of the results.
The structure constants in the classical limit have been compared 
with the classical structure constants and they agree. We have also 
verified the symmetry 
relations in section \ref{qlsec} for the quantum structure constants.
The relations derived in this section provide further checks. 
For quantum Lie algebras of
low rank we have checked that these relations are satisfied by our 
results. We have made the Mathematica notebooks containing these
calculations available on the Internet at 
http://www.mth.kcl.ac.uk/$\sim$delius/q-lie/

In the 
following we use the intertwiners  from adjoint $\otimes$ adjoint 
$\rightarrow$ adjoint and singlet $\rightarrow$ adjoint $\otimes$ adjoint. 
This second intertwiner is the inverse of the Killing 
form. Adopting the same notation as in section 
\ref{sect:kill} we denote this intertwiner by $B^{pq}$ and define it
so that $B_{pr}B^{rq}=\delta^q_p$. Thus in particular 
$B^{\alpha,-\alpha}=q^{-\rho\cdot\alpha}$. $B^{ij}$
(given in eqs. \eqref{killinvf}-\eqref{killinvl})
can be used to raise the indices which are lowered with $B_{ij}$.

We come now to the derivation of two new relations. They are obtained by 
constructing intertwiners singlet $\rightarrow $ adjoint and adjoint 
$\rightarrow$ adjoint. Further relations can easily be derived by the
same method.
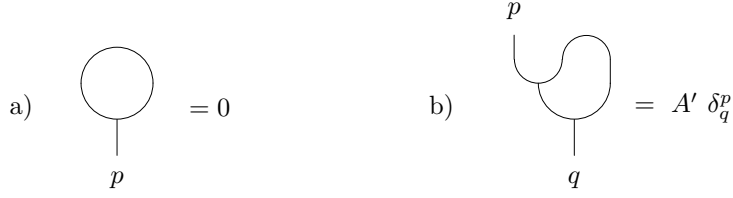
\begin{figure} 
\begin{center}
% Drawing generated by LaTeX-CAD 1.7 - requires latexcad.sty 
% (c) 1996 John Leis leis@usq.edu.au 
\begin{picture}(130,40)
\thinlines
\drawarc{20.0}{22.0}{12.0}{3.14}{6.28}
\drawarc{20.0}{22.0}{12.0}{6.28}{3.14}
\drawarc{98.0}{26.0}{8.0}{3.14}{6.28}
\drawarc{96.0}{22.0}{12.0}{6.28}{3.14}
\drawpath{102.0}{26.0}{102.0}{22.0}
\drawarc{90.0}{26.0}{8.0}{6.28}{3.14}
\drawpath{86.0}{30.0}{86.0}{26.0}
\drawpath{20.0}{16.0}{20.0}{10.0}
\drawpath{96.0}{16.0}{96.0}{10.0}
\drawcenteredtext{20.0}{6.0}{$p$}
\drawcenteredtext{96.0}{6.0}{$q$}
\drawcenteredtext{86.0}{34.0}{$p$}
\drawlefttext{32.0}{18.0}{$=0$}
\drawrighttext{6.0}{18.0}{a)}
\drawrighttext{76.0}{18.0}{b)}
\drawlefttext{106.0}{18.0}{$=~A'~ \delta_q^p$}
\end{picture}
\end{center}
\caption{Graphical representation of the relations. a) There is no 
non-zero intertwiner from singlet to adjoint. b) An intertwiner 
from adjoint to adjoint is proportional to the identity.
\label{symmfig}}
\end{figure}

\subsection{The intertwiner singlet $\rightarrow $ adjoint}

An intertwiner from the singlet to the adjoint is given in figure 
\ref{symmfig} a). This intertwiner should be zero and so we immediately 
arrive at a relation between the structure constants. Setting $p=k$ 
we have
\begin{equation}\label{prel1}
\sum_{\alpha}
-g_{\alpha}{}^{k} B^{\alpha,-\alpha}  +  f_{ij}^{~~k} B^{ij} = 0.
\end{equation}
We multiply this equation by $B_{km}$ and sum over $k$. We use
$B_{\alpha,-\alpha}=B^{\alpha,-\alpha}=q^{-\rho\cdot\alpha}$
and the relation \eqref{qant5}
between the $g_\alpha{}^k$ and the quantum roots to obtain
\begin{equation}\label{rel1}
\sum_{\alpha}l_{\alpha}(H_m) q^{-2 \rho \cdot \alpha} + 
\sum_i f_{mi}{}^i = 0.
\end{equation}

\subsection{The intertwiner adjoint $\rightarrow$ adjoint}

The intertwiner from adjoint to adjoint is given in figure 
\ref{symmfig} b). This intertwiner should be proportional to 
the identity map. Setting $p=q=\alpha$ we have
\begin{equation}
\sum_{\beta}N_{\alpha,\beta}N_{\alpha+\beta,-\beta}B^{\beta,-\beta}+
r_{\alpha}(H_i)r_{\alpha}(H_j) B^{ij}-
g_{\alpha}{}^{i} l_{\alpha}(H_i) B^{-\alpha,\alpha} =A'.
\end{equation}
Using once again the relation \eqref{qant5} between the $g_\alpha{}^i$
and the quantum roots, the relations between the $N_{\alpha\beta}$,
and the value of $B^{\alpha,-\alpha}$
the above equations become
\begin{equation}\label{rel2}
\sum_\beta
\left(N_{\alpha,\beta}\right)^2 q^{-2\rho\cdot\beta}+
r_{\alpha}(H_i)r_{\alpha}(H_j) B^{ij}+ 
l_{\alpha}(H_i)l_{\alpha}(H_j) B^{ij} =A',\ \ \ \forall \alpha\in R.
\end{equation}
 
Setting $p=i$ and $q=j$ in Figure \ref{symmfig} b) we obtain
\begin{equation}
-\sum_\alpha\ l_\alpha(H_i)\,B^{\alpha,-\alpha}\,g_\alpha{}^j+
f_{ik}{}^m\,B^{kl}\,f_{ml}{}^j\,=\,A'~\delta_i^j.
\end{equation}
which can be reexpressed as
\begin{equation}\label{rel3}
\sum_\alpha\ l_{\alpha}(H_i)l_{\alpha}(H_j)+
f_i{}^{ml}\,f_{mlj}\,=\,A'~B_{ij}.
\end{equation}

\section{Discussion\label{sect:disc}}
We have shown how to calculate the structure constants of the 
quantum Lie algebras associated to 
$B_l$, $C_l$ and $D_l$. These calculations were rendered manageable 
by the observation that the quantum structure constants are just the 
inverse Clebsch-Gordan coefficients for adjoint $\otimes$ adjoint 
$\rightarrow$ adjoint. The structure constants satisfy the symmetries 
discovered in \cite{Del96}. We have introduced an ad-invariant 
Killing form and shown that it is proportional to the intertwiner 
from adjoint 
$\otimes$ adjoint $\rightarrow$ singlet. Because the composition of 
intertwiners is also an intertwiner we were able to calculate many 
intertwiners indirectly. For example we calculated the Killing form 
by building the intertwiner from adjoint $\otimes$ adjoint 
$\rightarrow$ singlet in terms of the intertwiners from vector 
$\otimes$ vector 
into adjoint and singlet. This meant that we didn't have to evaluate 
the usual trace over the adjoint representation.

As is well known, the structure constants of the simple complex Lie 
algebras are determined entirely in terms of their simple roots.
Eventually we would hope to arrive at a similar result for the
quantum Lie algebras. In this paper we have come one step closer
to this goal by our observation that the 
structure constants $f_{ijk}$ of the Cartan subalgebra are completely
determined in terms of the quantum roots according to eq. \eqref{fr}.
The Killing form $B_{ij}$ is expressed in terms of the $l_\a$ by
eq. \eqref{br}.
Already in \cite{Del96} it was found that the left quantum roots
$l_\a$ for positive $\a$ are enough to determine those for negative
$\a$ by eq. \eqref{qant2}, the $r_\a$ by \eqref{qant1} and the $g_\a{}^k$
by eq. \eqref{qant5}. Thus now all quantum Lie bracket relations are
determined by the left quantum roots $l_\a$ for positive
$\a$ and the $N_{\a\beta}$. What is still missing is a
deeper understanding
of the $N_{\a\beta}$ and of how to obtain the higher quantum roots from
the simple quantum roots. For recent progress see \cite{Gar98}.

The expressions for the quantum roots in eqs. \eqref{ba}, \eqref{ca}
and \eqref{da} are unexpectedly simple. If one writes the
classical expressions for the roots in the same form, one notices
that generically the quantum expressions are obtained from these 
by replacing every $1$ by a power of $q$ and every $2$ by $(q+q^{-1})$ 
times a power of $q$. Thus in particular $l_\a(H_i)\neq 0$ if and only
if classically $\a(H_i)\neq 0$. There is however one exception to
this simplicity: in $(B_l)_h$ we have found that $l_{\epsilon_k}(H_l)
\neq 0$ also for $k<l-1$.

Also the matrices \eqref{kill1}-\eqref{kill3}
describing the quantum Killing form on the Cartan subalgebras
are surprisingly simple. In particular we find that only those
entries in the matrices are non-zero which are also non-zero classically.

%\begin{eqnarray}\label{new}
%{}[H_i,X_\alpha]_h&=&l_\alpha(H_i)\,X_{\alpha},\ \ \ \ \ \ 
%{}[X_\alpha,H_i]_h=-\widetilde{l_{\alpha}}(H_i)\,X_{\alpha},\nonumber\\
%{}[H_i,H_j]_h&=&B_{ij}\left(l^a_{\a_j}(H_i)\,H^i+
%l^a_{\a_i}(H_j)\,H^j\right),\ \ \ ~i\neq j,\nonumber\\
%{}[X_\alpha,X_{-\alpha}]_h&=&\sum_k\,B_{\a,-\a}\,l_\alpha(H_k)\,H^k,
%\nonumber\\
%{}[H_i,H_i]_h&=&\sum_k\,B_{ik}\,l^a_{\a_k}(H_i)\,H^k,\nonumber\\
%\bigl[X_{\alpha} , X_{\beta} \bigr]_h & = & \left\{ 
%\begin{array}{ll}
%N_{\alpha,\beta} ~X_{\alpha+\beta} & 
%\mbox{    for  } \alpha+\beta ~\in~ R \\
% 0 & \mbox{     otherwise}. \end{array} \right. 
%\end{eqnarray}
% 
There still remain a lot of unanswered questions. In particular:
What is a good axiomatic setting for the theory of
quantum Lie algebras. How should one $q$-deform the 
Jacobi identity? What characterizes the quantum root system? What 
are $q$-Weyl reflections? How does one define representations of
these non-associative algebras? And many more.
For a more complete bibliography and more recent results on
quantum Lie algebras see
http://www.mth.kcl.ac.uk/$\sim$delius/q-lie/

\appendix
\section{Quantized enveloping algebras}
\label{hopf}
For an introduction to Lie algebras consult \cite{Sam69,Hum70} 
and for quantized enveloping algebras consult \cite{Cha94}. 
\begin{definition}
Let $\mathfrak{g}$ be a finite-dimensional simple complex Lie 
algebra with symmetrizable Cartan matrix $a_{ij}$. The 
\dem{quantized enveloping algebra} $U_h(\mathfrak{g})$ is the 
unital associative algebra over $\mathbb C[[h]]$ (completed in 
the $h$-adic topology) with generators $x_i^+,\ x_i^-,\ h_i$, 
$1 \le i \le \text{rank}(\mathfrak{g})$ and relations
\begin{gather} \label{defrel}
 h_i h_j = h_j h_i,\ \ \ ~~
h_i x_j^\pm - x_j^\pm h_i = \pm a_{ij} x_j^\pm , \nonumber\\
x_i^+ x_j^- -x_j^- x_i^+ =  \delta_{ij}
\ \frac{q_i^{h_i} - q_i^{-h_i}}{q_i- q_i^{-1}},\\
\sum_{k=0}^{1-a_{ij}} (-1)^k
\left[ \begin{array}{c} 1-a_{ij} \\ k \end{array} \right]_{q_i}
(x_i^\pm)^k x_j^\pm (x_i^\pm)^{1-a_{ij}-k} = 0 \qquad i \not= j.
\nonumber
\end{gather}
Here $\left[\begin{array}{c}a\\b\end{array}\right]_q$ are the 
$q$-binomial coefficients. 

We have defined $q_i = e^{d_i h}$ where the $d_i$ are chosen so 
that $d_i a_{ij}$ is a symmetric matrix. We choose 
$d_i= \alpha_i^2/2$ where the simple roots are as given at the 
beginning of section \ref{sect:cgc}. An alternative convention 
is to choose the $d_i$ to be coprime integers. In the case of 
the algebra $B_l$ these two conventions differ and our conventions 
lead to $d_i=1$ for $i=1,~..,~l-1$ and $d_l=\frac{1}{2}$. The 
Cartan matrix is defined to be $a_{ij}=2 \alpha_i \cdot 
\alpha_j/\alpha_i^2$.
\end{definition}
The Hopf algebra structure of $U_h(\mathfrak{g})$ is given by the 
comultiplication $\Delta:U_h(\mathfrak{g}) \rightarrow 
U_h(\mathfrak{g}) \hat{\otimes} U_h(\mathfrak{g})$ 
defined by
\begin{align} 
\Delta(h_i) &= h_i \hat{\otimes} 1 + 1 \hat{\otimes} h_i, \\
\Delta(x_i^\pm) &= x_i^\pm \hat{\otimes} q_i^{-h_i/2} +
 q_i^{h_i/2} \hat{\otimes} x_i^\pm,
\end{align}
and the antipode $S$ and counit $\epsilon$ defined by
\begin{equation}
S(h_i)= - h_i, \ \ \ 
S(x_i^\pm) = - q_i^{\mp 1}\,x_i^\pm ,\ \ \ 
\epsilon(h_i) = \epsilon(x_i^\pm) = 0.
\end{equation}
\begin{definition} The Cartan involution $\theta:~U_h(\mathfrak{g}) 
\rightarrow U_h(\mathfrak{g})$ is given by the formulae
\begin{equation}
\theta(x_i^{\pm})~=~x_i^{\pm},~\theta(h_i)~=~-h_i.
\end{equation}
\end{definition}
It is an Hopf-algebra isomorphism: $U_h(\mathfrak{g}) \rightarrow 
U^{\text{op}}_h(\mathfrak{g})$ where $U^{\text{op}}_h(\mathfrak{g})$ is the
opposite Hopf algebra, whose Hopf structure is described by the 
opposite coproduct 
$\Delta^{\text{op}}$ and the inverse antipode $S^{-1}$.
\begin{definition}
$q$-conjugation
\mbox{$\sim: \mathbb C[[h]] \rightarrow \mathbb C[[h]]$}, 
$a\mapsto\tilde{a}$ is the $\mathbb C$-linear ring
automorphism defined by $\tilde{h}=-h$.
\end{definition}

\section{Modified Schur's lemma}
\label{schur}
\begin{lemma}[Schur's lemma]
Let $V[[h]]$ and $W[[h]]$ be two finite-di\-men\-sio\-nal indecomposable
$U_h(\mathfrak{g})$-modules and let $f,g:V[[h]]\rightarrow W[[h]]$
be two $U_h(\mathfrak{g})$-module homomorphism. Then 
\begin{enumerate}
\item if $f~(mod~h) \neq 0$ then $f$ is an isomorphism.
\item $\exists c \in \mathbb C[[h]]$ such that $f=cg$ or $g=cf$.
\end{enumerate}
\end{lemma}
\begin{proof}
(1) Ker$(f)$ is a submodule of $V[[h]]$. But $V[[h]]$ is an 
indecomposable $U_h(\mathfrak{g})$-module and so Ker$(f)$ must 
be of the form $c~V[[h]]$ for some non-invertible $c \in 
\mathbb C[[h]]$. However if Ker$(f)$ has this form then 
$f(c~x)=c~f(x)=0$ $\forall x \in V[[h]]$, i.e. Ker$(f) = V[[h]]$. 
Therefore Ker$(f)=0$. Equally Im$(f)$ is a submodule of $W[[h]]$ 
and so 
(2) Let $v_0$ and $w_0$ be the highest weight states in $V$ and 
$W$. Then because $f,g$ are $U_h(\mathfrak{g})$-homomorphisms 
$f(v_0)$ and $g(v_0)$ must also be highest weight states in 
$W[[h]]$, that is $\exists ~c_1,c_2 \in \mathbb C[[h]]$ such 
that $f(v_0)=c_2 w_0$ and $g(v_0)=c_1 w_0$. Then 
$(c_1 f - c_2 g)(v_0)=0$ which means that $c_1 f - c_2 g$ is not 
an isomorphism and so by the first part of the lemma 
$c_1 f - c_2 g = 0~(mod~h)$. By the same argument 
$h^{-1}(c_1 f - c_2 g)=0$, etc... Hence $c_1 f = c_2 g$.
\end{proof}

\section{Clebsch-Gordan coefficients}
\label{cleb}
Let $(\pi^{\mu},V^{\mu})$ and $(\pi^{\nu},V^{\nu})$ be two indecomposable 
$U_h(\mathfrak{g})$-modules. Consider an 
indecomposable $U_h(\mathfrak{g})$-module $(\pi^{\lambda},V^{\lambda})$ 
homomorphically embedded in $V^{\mu} \otimes 
V^{\nu}$. As a basis for $V^{\lambda}$ we choose \{ $v^a_{\lambda}$ \}. 
So we have for the action of $U_h(\mathfrak{g})$ on $V^{\lambda}$
\begin{equation}
\pi^{\lambda}(x) v^c_{\lambda}~~=~~ v^d_{\lambda} 
\pi^{\lambda}_d{}^c(x) ~=~  \pi^{\lambda}_d{}^c(x)~
C^{\mu \nu}_{a' b'}|^d_{\lambda}~ v^{a'}_{\mu} \otimes v^{b'}_{\nu}
\end{equation}
where \{ $v^{a'}_{\mu} \otimes v^{b'}_{\nu}$ \} is the natural basis 
on $V^{\mu} \otimes V^{\nu}$ and $C^{\mu \nu}_{a'b'}|^d_{\lambda}$ 
are the Clebsch-Gordan coefficients describing the embedding 
$V^{\lambda} \rightarrow V^{\mu} \otimes V^{\nu}$. The action of 
$U_h(\mathfrak{g})$ on $V^{\mu} \otimes V^{\nu}$ is defined using 
the coproduct, i.e.  
\begin{eqnarray}
\pi^\lambda(x) ~ v^c_{\lambda}  & = &  
C^{\mu \nu}_{ab} |^c_{\lambda}~
(\pi^{\mu} \otimes \pi^{\nu})(\Delta(x))~
 (v^{a}_{\mu} \otimes v^{b}_{\nu})\\ \nonumber
 & = & C^{\mu \nu}_{ab} |^c_{\lambda}~ (v^{a'}_{\mu} \otimes 
v^{b'}_{\nu})~ (\pi^{\mu a}_{a'} \otimes \pi^{\nu b}_{b'})(\Delta(x)) 
\end{eqnarray}
These two actions of $U_h(\mathfrak{g})$  coincide and so the 
Clebsch-Gordan coefficients satisfy the intertwiner property 
\begin{equation}
\pi^{\lambda}(x)^{~ c }_d~ C^{\mu \nu}_{a'b'} |^d_{\lambda} ~~=~~ 
C^{\mu \nu}_{ab}|_{\lambda}^c ~(\pi^{\mu a}_{a'} \otimes 
\pi^{\nu b}_{b'})(\Delta(x)).
\end{equation}

%\section{Adjoint representation matrices}
%
%Here we give the adjoint representation matrices  of the $U_h(\lie)$
%generators $\{x^+_i,x^-_i,h_i\}_{i=1,..,l}$. There are two 
%important points to comment 
%on. Firstly the representation matrices of the Cartan elements 
%are independent of $q$, only the representation matrices of the root 
%vectors are $q$-deformed, and secondly in our basis $\adj{x^+_i}\neq
%\adj{x^-_i}^\dagger$.
%
%\subsection{$\mathfrak{g}= B_l$}
%\begin{eqnarray}
%\adj{h_i}X_\alpha&=&\frac{2\alpha_i\cdot\alpha}{\alpha_i^2}\ X_\alpha,
%\nonumber\\
%\adj{h_i}H_j&=&0,\nonumber\\
%\adj{x^+_i}X_{\alpha_{jk}} 
%& = & 
%\end{eqnarray}
%

\section{The $(C_2)_h$ algebra}
In this appendix we compare our results for $(C_2)_h$ with the
results given in \cite{Del96}. The Cartan matrix for $C_2$ is
\begin{gather}
\begin{pmatrix} 2 & -2 \\ -1 & 2 \end{pmatrix} 
\end{gather}
and the positive roots are $\alpha_1,~\alpha_2,~\alpha_1+\alpha_2$, and 
$2\alpha_1+\alpha_2$. The left quantum roots are given by
\begin{eqnarray}
l_{\alpha_1}(H_1,H_2) & = &  ((q+q^{-1})(q^2-1+q^{-2})q,-q^{-3})),
\nonumber \\
l_{\alpha_2}(H_1,H_2) & = &  (-(q+q^{-1})q,q^2(q+q^{-1})), \nonumber \\
l_{\alpha_1+\alpha_2}(H_1,H_2) & = &  ((q+q^{-1})(q-q^{-1})q^2,
q)), \nonumber \\
l_{2\alpha_1+\alpha_2}(H_1,H_2) & = &  ((q+q^{-1})q^3,0).
\end{eqnarray}
We agree with the previous results in \cite{Del96} if we first 
$q$-conjugate our results and then make the following transformations:
\begin{eqnarray}
H_1 & \longrightarrow & l\, \frac{(q^{-2}+q^{2}-1)}{q+q^{-1}}h_1, 
\nonumber \\
H_2 & \longrightarrow & l\, h_2, 
\end{eqnarray}
where $l$ 
is as defined in \cite{Del96} and $h_1,~ h_2$ are the Cartan subalgebra 
generators used in \cite{Del96}. 

In the table below we give 
$N_{\alpha,\beta}/(q+q^{-1})^{1/2}$ for $\alpha$ positive. 
The rows are labelled by $\alpha$ and the columns by $\beta$. 
The $N_{-\alpha,\beta}$ are equal to $-\tilde{N}_{\alpha,-\beta}$.
\[
\begin{array}{r|cccccccc}
&\begin{array}{c}2\a_1\\~+\a_2\end{array}\!\!&
\!\!\begin{array}{c}\a_1\,\\~+\a_2\end{array}&
{}~~\a_2&~~\a_1&~-\a_1&~-\a_2&
\begin{array}{c}-\a_1\\~-\a_2\end{array}&
\!\!\begin{array}{c}-2\a_1\\~-\a_2\end{array}\!\!\\
\hline
2\alpha_1 + \alpha_2 &  0 & 0 & 0 & 0 & q^{-2} & 0 & -q^{-2}  & 0 \\
\alpha_1 + \alpha_2 & 0 & 0 & 0 & q^{-1} & q^{-3} & -1 & 0 & 
-q^{-2} \\
\alpha_2 & 0 & 0 & 0 & q^{-2} & 0 & 0 & -1 & 0 \\
\alpha_1 & 0 & -q & -q^2  & 0 & 0 & 0 & q^{-3} & q^{-2} \\
\end{array}
\]

We get agreement with \cite{Del96} if we first $q$-conjugate 
and then make the following transformations:
\begin{eqnarray}
X_{\pm \alpha_1} & \longrightarrow \mp  \xi\, X_{\pm \alpha_1}, \nonumber \\
X_{\pm \alpha_2} & \longrightarrow \mp \xi\, X_{\pm \alpha_2}, \nonumber \\
X_{\pm (\alpha_1+\alpha_2)} & \longrightarrow \mp \xi\, X_{\pm 
(\alpha_1+\alpha_2)}, \nonumber \\
X_{\pm (2 \alpha_1+\alpha_2)} & \longrightarrow \pm \xi\, X_{\pm 
(2 \alpha_1+\alpha_2)},
\end{eqnarray}
where $\xi=-(q+q^{-1})^{1/2}(q^2-1+q^{-2})$. The algebra $(C_2)_h$ is 
isomorphic to the algebra $(B_2)_{2h}$. The change in $h$ is due to our 
choice of conventions for $d_i$ in the definition of $U_h(\mathfrak{g})$ 
in Appendix \ref{hopf}.


\begin{thebibliography}{99}

\bibitem {Asc93} P. Aschieri, L. Castellani, {\it An introduction to
noncommutative differential geometry on quantum groups},
Int. J. Mod. Phys. {\bf A8} (1993) 1667.

\bibitem {Ber90} D. Bernard, {\it Quantum Lie Algebras and
Differential Calculus on Quantum Groups}, Prog. Theo. Phys. Suppl. {\bf 102}
(1990) 49.

\bibitem {Cha94}        V. Chari, A. Pressley,   {\it A Guide to Quantum
Groups}, Cambridge University Press {\bf } (1994).

\bibitem {Del96} G.W. Delius, A. H\"uffmann {\it On quantum Lie algebras 
and quantum root systems}, q-alg/9506017,
J. Phys {\bf A 29} (1996) 1703.

\bibitem{Del95} G.W. Delius, A. H\"uffmann, M.D. Gould, Y.-Z. Zhang,
{\it The quantum Lie algebras associated to $U_h(gl_n)$ and
$U_h(sl_n)$},  q-alg/9508013, J. Phys. {\bf A 29} (1996) 5611.

\bibitem {Del97} G.W. Delius 
and M. Gould, {\it Quantum Lie algebras, their existence, uniqueness and 
$q$-antisymmetry}, q-alg/9605025, Commun. Math. Phys {\bf 185} (1997) 709.

\bibitem {DelP96} G.W. Delius, {\it The problem of differential calculus 
on quantum groups}, q-alg/9608010,
Proc. Quantum groups and Integral Systems, Prague 1996.

\bibitem {Dri85} V.G. Drinfel'd, {\it Hopf algebras and the quantum
Yang-Baxter equation}, Sov. Math. Dokl. {\bf 32} (1985) 254, {\it Quantum 
Groups}, Proc. Int. Congr. Math., Berkeley {\bf } (1986) 798.

\bibitem {Dri90} V.G. Drinfel'd, {\it On almost cocommutative Hopf 
algebras}, Leningrad Math. J., Vol. 1 {\bf No. 2} (1990). 

\bibitem {Fad87} L.D. Faddeev, P.N. Pyatov, {\it The differential calculus 
on quantum linear groups}, hep-th/9402070, 
AMS Transl. (2) Vol. 175 (1996) 35-47; {\it 
Quantization of differential calculus
on linear groups}, in 'Problems in Modern Theoretical Physics',
Ed. A.P.Isaev, JINR Publ. Dept. 96-212, Dubna, 1996, pp. 19-43.

\bibitem{Gar98} C. Gardner, preprint in preparation, see 
http://www.mth.kcl.ac.uk/$\sim$delius/q-lie/

\bibitem {Hum70} J. E. Humphreys, {\it Introduction to Lie Algebras and 
Representation Theory}, Springer (1970).

\bibitem {Ros90} M. Rosso,       {\it Analogues de la forme de Killing
et du theoreme  d'Harish-Chandra pour les groupes quantiques},
Ann. scient. Ec. Norm. Sup. {\bf 23} (1990) 445.

\bibitem {Sam69} H. Samelson,    {\it Notes on Lie Algebras}, Van
Nostrand {\bf } (1969).

\bibitem {Sud95} V. Lyubashenko, A. Sudbery, {\it Quantum Lie Algebras of 
Type} $A_n$, q-alg/9510004.

\bibitem {Swe69} M.E. Sweedler, {\it Hopf algebras}, Benjamin, New 
York (1996).

\bibitem {Wol91} S. Wolfram, {\it Mathematica}, Addison-Wesley Publishing 
Co. (2nd ed.) (1991).

\bibitem{Wor89} S.L. Woronowicz, {\it Differential Calculus on Compact
Matrix Pseudogroups (Quantum Groups)}, Commun. Math. Phys. {\bf 122}
(1989) 125.



\end{thebibliography}
\end{document}